\documentclass[11pt]{article}

\usepackage{amscd,amsmath, amssymb, fancyhdr, epsfig,url,color}

\numberwithin{equation}{section}

\newcommand{\version}{v. 1.2.1, August 18, 2017}

\def\eqref#1{(\ref{#1})}
\newcommand{\goth}{\mathfrak}

\newcommand{\arrow}{{\:\longrightarrow\:}}
\newcommand{\Z}{{\mathbb Z}}
\newcommand{\C}{{\mathbb C}}

\newcommand{\R}{{\mathbb R}}
\newcommand{\Q}{{\mathbb Q}}

\def\1{\sqrt{-1}\:}

\newcommand{\cntrct}                % contraction with a vector field
{\hspace{2pt}\raisebox{1pt}{\text{$\lrcorner$}}\hspace{2pt}}

\makeatletter
\def\x@arrow{\DOTSB\Relbar}
\def\xlongequalsignfill@{\arrowfill@\x@arrow\Relbar\x@arrow}
\newcommand{\xlongequal}[2][]{%
        \ext@arrow 0099\xlongequalsignfill@{#1}{#2}}
\def\xlongrightarrowfill@{\arrowfill@\relbar\relbar\longrightarrow}
\newcommand{\xlongrightarrow}[2][]{%
        \ext@arrow 0099\xlongrightarrowfill@{#1}{#2}}
\makeatother

\newcommand{\Vol}{\operatorname{Vol}}

\newcommand{\Lie}{\operatorname{Lie}}
\newcommand{\Sp}{\operatorname{Sp}}

\newcommand{\Kah}{\operatorname{\sf Kah}}

\newcounter{Mycounter}[section]
\newcounter{lemma}[section]
\renewcommand{\thelemma}{{Lemma \thesection.\arabic{lemma}}}
\newcommand{\lemma}{%
    \setcounter{lemma}{\value{Mycounter}}
    \refstepcounter{lemma}
    \stepcounter{Mycounter}
    {\noindent \bf \thelemma:\ }}

\newcounter{claim}[section]
\newcounter{sublemma}[section]
\newcounter{corollary}[section]
\renewcommand{\thecorollary}{{Corollary \thesection.\arabic{corollary}}}
\newcommand{\corollary}{%
    \setcounter{corollary}{\value{Mycounter}}
    \refstepcounter{corollary}
    \stepcounter{Mycounter}
    {\noindent \bf \thecorollary:\ }}

\newcounter{theorem}[section]
\renewcommand{\thetheorem}{{Theorem \thesection.\arabic{theorem}}}
\newcommand{\theorem}{%
    \setcounter{theorem}{\value{Mycounter}}
    \refstepcounter{theorem}
    \stepcounter{Mycounter}
    {\noindent \bf \thetheorem:\ }}

\newcounter{conjecture}[section]
\newcounter{proposition}[section]
\renewcommand{\theproposition}
      {{Proposition \thesection.\arabic{proposition}}}
\newcommand{\proposition}{%
    \setcounter{proposition}{\value{Mycounter}}
    \refstepcounter{proposition}
    \stepcounter{Mycounter}
    {\noindent \bf \theproposition:\ }}

\newcounter{definition}[section]
\renewcommand{\thedefinition}
      {{Definition~\thesection.\arabic{definition}}}
\newcommand{\definition}{%
    \setcounter{definition}{\value{Mycounter}}
    \refstepcounter{definition}
    \stepcounter{Mycounter}
    {\noindent \bf \thedefinition:\ }}

\newcounter{example}[section]
\renewcommand{\theexample}{{Example \thesection.\arabic{example}}}
\newcommand{\example}{%
    \setcounter{example}{\value{Mycounter}}
    \refstepcounter{example}
    \stepcounter{Mycounter}
    {\noindent \bf \theexample:\ }}

\newcounter{remark}[section]
\renewcommand{\theremark}{{Remark \thesection.\arabic{remark}}}
\newcommand{\remark}{%
    \setcounter{remark}{\value{Mycounter}}
    \refstepcounter{remark}
    \stepcounter{Mycounter}
    {\noindent \bf \theremark:\ }}

\newcounter{problem}[section]
\newcounter{question}[section]
\makeatletter

\@addtoreset{equation}{section} \@addtoreset{footnote}{section} \makeatother

\def\blacksquare{\hbox{\vrule width 5pt height 5pt depth 0pt}}
\def\endproof{\blacksquare}

\newcommand{\N}{{\mathbb{N}}}

\newcommand{\ba}{{\bar{a}}}

\newcommand{\bx}{{\bar{x}}}

\newcommand{\cE}{{\mathcal{E}}}

\newcommand{\cK}{{\mathcal{K}}}

\newcommand{\cS}{{\mathcal{S}}}

\newcommand{\hI}{{\hat{I}}}
\newcommand{\hJ}{{\hat{J}}}

\newcommand{\hN}{{\hat{N}}}

\newcommand{\heta}{{\hat{\eta}}}

\newcommand{\hpi}{{\hat{\pi}}}

\newcommand{\goh}{{\mathfrak{h}}}

\newcommand{\Span}{{{\rm Span}\,}}

\newcommand{\tI}{{\widetilde{I}}}
\newcommand{\tJ}{{\widetilde{J}}}

\newcommand{\tM}{{\widetilde{M}}}

\newcommand{\tU}{{\widetilde{U}}}

\newcommand{\tomega}{{\widetilde{\omega}}}

\newcommand{\tphi}{{\widetilde{\phi}}}

\newcommand{\talpha}{{\widetilde{\alpha}}}

\newcommand{\tOmega}{{\widetilde{\Omega}}}

\begin{document}

%\color{black}
%%%%%%%%%%%%%%%%%%%%%%%%%%%%%%%%%%%%%%%%%%%%%%%%%%%%%%%%%%%%
\begin{center}
{\LARGE\bf
Unobstructed symplectic packing by ellipsoids
for tori and hyperk\"ahler manifolds\\
\ \\
}
%%%%%%%%%%%%%%%%%%%%%%%%%%%%%%%%%%%%%%%%%%%%%%%%%%%%%%%%%%%%

Michael Entov\footnote{Partially supported by the Israel Science
Foundation grant $\#$ 1096/14.}, Misha Verbitsky\footnote{
  Partially supported by the Russian Academic Excellence Project '5-100'.
 }

\end{center}

%%%%%%%%%%%%%%%%%%%%%%%%%%%%%%%%%%%%%%%%%%%%%%%%

\hfill
\hfill
%%%%%%%%%%%%%%%%%%%%%%%%%%%%%%%%%%%%%%%%%%%%%%%%
%{\small \hspace{0.10\linewidth}
%\begin{minipage}[t]{0.85\linewidth}
%{\bf Abstract} \\

\begin{abstract}
\small Let
$M$ be a closed symplectic manifold of volume $V$. We say that the symplectic packings of $M$
by ellipsoids are unobstructed if any collection
of disjoint symplectic ellipsoids (possibly of different sizes) of total volume
less than $V$ admits a symplectic embedding to $M$.
We show that the symplectic packings by ellipsoids are unobstructed for all even-dimensional tori equipped with K\"ahler symplectic forms
and all closed hyperk\"ahler manifolds of maximal holonomy, or, more generally, for closed Campana simple manifolds (that is, K\"ahler manifolds that are not unions of their complex subvarieties),
as well as for any closed K\"ahler manifold which is a
limit of Campana simple manifolds in a smooth deformation.
The proof involves the construction of a K\"ahler resolution of a K\"ahler orbifold with isolated singularities
and relies on the results of Demailly-Paun and Miyaoka on K\"ahler cohomology classes.

\end{abstract}

%\end{minipage}
%}
%%%%%%%%%%%%%%%%%%%%%%%%%%%%%%%%%%%%%%%%%%%%%%%%

{\small
\tableofcontents
}

%%%%%%%%%%%%%%%%%%%%%%%%%%%%%%%%%%%%%%%%%%%%%%%%

\section{Introduction}

%%%%%%%%%%%%%%%%%%%%%%%%%%%%%%%%%%%%%%%%%%%%%%%%

The symplectic packing problem is one of the
central
problems of symplectic topology -- it concerns the existence of symplectic embeddings of
a union of disjoint copies of various (possibly different) sizes of a particular standard shape (ball, ellipsoid, polydisk etc.) into a given $2n$-dimensional symplectic manifold $(M,\omega)$.
An immediate obstruction to such symplectic embeddings is given by the symplectic volume. It has been known since the pioneering work by Gromov \cite{_Gromov_} that there might be additional obstructions
coming from pseudo-holomorphic curves in $(M,\omega)$ -- a symplectic rigidity phenomenon. In \cite{_EV-balls_} we prove an opposite, symplectic flexibility, claim for the symplectic packing of {\it K\"ahler} manifolds by balls: in the case when $(M,\omega)$ is a K\"ahler manifold admitting ``really few" (genuinely, not pseudo-) holomorphic subvarieties (such a K\"ahler manifold is called Campana simple), or if it can be approximated by Campana simple manifolds, then there are no obstructions to symplectic embeddings of disjoint unions of balls into $(M,\omega)$ apart from the volume. In this paper we extend this flexibility result to symplectic packings of K\"ahler manifolds by ellipsoids.

Let us say a few words about the method of the proof.
In the case of the symplectic packings by balls,  McDuff and Polterovich \cite{_McD-Polt_} reduced the
question about symplectic embeddings of unions of $k$ balls into a symplectic
manifold $(M,\omega)$ to a question about
the structure of the symplectic cone in the cohomology of
a blow-up of $M$ at $k$ points. In the same paper they
showed that symplectic packings of K\"ahler manifolds by balls
are deeply related to algebraic geometry that
allows sometimes to describe the shape of the K\"ahler cone in the cohomology of a
K\"ahler manifold.
In \cite{_EV-balls_}
we proved the above-mentioned flexibility result for the symplectic packings by balls using
the results of McDuff-Polterovich along
with several strong results from complex geometry -- in particular,
the Demailly-Paun theorem
\cite{_Demailly_Paun_}
describing completely the K\"ahler cone of a closed K\"ahler manifold.
(A similar approach was previously used by Latschev-McDuff-Schlenk
\cite{_LMcDS_} in the case when $(M,\omega)$ is a K\"ahler torus of real dimension 4).

The problem with extending this kind of argument to the study of symplectic
packing by ellipsoids is that instead of the usual blow-ups of a symplectic manifold
one has to consider weighted blow-ups which, unlike the usual blow-ups, produce
not smooth manifolds but orbifolds, where the results and the techniques
used to describe the symplectic/K\"ahler cone of the usual blow-up do not apply.
(An orbifold version of the Demailly-Paun theorem seems to be true but its proof is not published, as far as we can ascertain).
Therefore we use an indirect argument where
the Demailly-Paun theorem is applied not to the orbifold but to its K\"ahler resolution which is a smooth manifold.

Let us note that McDuff \cite{_McDuff-ellipsoids_} invented an approach to the study of symplectic packing of balls by ellipsoids in real dimension 4 which avoids dealing with orbifolds and which is based on an ingenious trick reducing the problem about symplectic embeddings of an ellipsoid
to the problem about symplectic embeddings of a certain disjoint union of balls.
F.Schlenk has independently proved the flexibility results for symplectic packings of tori and hyperk\"ahler manifolds by ellipsoids
{\it in dimension} 4 \cite{_Schlenk-survey_} by combining McDuff's method with our flexibility results for symplectic packing by balls \cite{_EV-balls_}. The proof that we give below works in all dimensions.

%%%%%%%%%%%%%%%%%%%%%%%%%%%%%%%%%%%%%%%%%%%%%%%%%%%%%%%%%%%%
\section{Main results}
\label{_USP_Section_}
%%%%%%%%%%%%%%%%%%%%%%%%%%%%%%%%%%%%%%%%%%%%%%%%%%%%%%%%%%%%

%%%%%%%%%%%%%%%%%%%%%%%%%%%%%%%%%%%%%%%%%%%%%%%%%%%%%%%%%%%%%%%%%%%%%%%%
\subsection{Preliminaries}
\label{_Preliminaries_Section_}
%%%%%%%%%%%%%%%%%%%%%%%%%%%%%%%%%%%%%%%%%%%%%%%%%%%%%%%%%%%%%%%%%%%%%%%%

\noindent {\bf Symplectic and complex structures}.
We view complex structures as tensors, that is, as integrable almost complex structures.

We say that an almost complex structure $J$ and a differential 2-form
$\omega$ on a smooth manifold $M$ are {\bf
compatible} with
each other if $\omega (\cdot, J\cdot)$ is a
$J$-invariant Riemannian metric on $M$.

The compatibility between a {\it complex} structure $J$ and a
symplectic form $\omega$ means exactly that $\omega (\cdot, J\cdot)
+ i\omega (\cdot, \cdot)$ is a K\"ahler metric on $M$.

We call
a symplectic form {\bf K\"ahler}, if it is compatible with {\it some}
complex structure.

A degree-2 real cohomology class of a complex manifold $(M,J)$ is called {\bf K\"ahler (with respect to $J$)} if it can be realized
by a K\"ahler form compatible with $J$. Such classes form an open cone that will be denoted by $\Kah (M,J)\subset H^2 (M;\R)$.

We will say that a complex structure is of {\bf K\"ahler type} if
it is compatible with {\it some} symplectic form.

\bigskip
\noindent {\bf Symplectic forms on tori}. Consider a torus $T^{2n} =
\R^{2n}/\Z^{2n}$ and let $\pi: \R^{2n}\to \R^{2n}/\Z^{2n}=T^{2n}$
be the natural projection.

%A differential form
%on
%$\R^{2n}$ is called {\bf linear}
%if it has constant coefficients with respect to the
%standard coordinates on $\R^{2n}$.
%A linear differential form
%on $\R^{2n}$
%descends under $\pi$ to a differential form
%on $T^{2n}$. We call a differential form
%on $T^{2n}$ {\bf linear} if it can be obtained in this way.
The K\"ahler forms on $T^{2n}$ are exactly the ones that can be mapped
by a diffeomorphism of $T^{2n}$ to a symplectic form whose lift by $\pi$ to $\R^{2n}$
has constant coefficients with respect to the
standard coordinates on $\R^{2n}$
(see e.g. \cite[Proposition 6.1]{_EV-balls_}).

\bigskip
\noindent {\bf Hyperk\"ahler manifolds}. There are several
equivalent definitions of a hyperk\"ahler manifold.
Since we study hyperk\"ahler manifolds from the symplectic
viewpoint, here is a definition which is close in spirit to
symplectic geometry: A {\bf hyperk\"ahler manifold} is a manifold
equipped with three complex structures $I_1,I_2,I_3$ satisfying the
quaternionic relations and three symplectic forms
$\omega_1,\omega_2,\omega_3$ compatible, respectively, with
$I_1,I_2,I_3$, so that the three Riemannian metrics $\omega_i
(\cdot, I_i\cdot)$, $i=1,2,3$, coincide. Such a collection of
complex structures and symplectic forms on a manifold is called a
{\bf hyperk\"ahler structure} and will be denoted by
$\goh=\{ I_1, I_2,I_3, \omega_1,\omega_2,\omega_3\}$.

We will say that a symplectic form is
{\bf hyperk\"ahler}
and a complex structure is of {\bf hyperk\"ahler type},
if they appear in {\it some} hyperk\"ahler
structure. In particular, any hyperk\"ahler symplectic form is
K\"ahler
and any complex structure of hyperk\"ahler type is also of K\"ahler type.

A hyperk\"ahler manifold $(M,\goh)$ is called {\bf irreducible
holomorphically symplectic (IHS)} if $\pi_1 (M)=0$ and $\dim_\C
H^{2,0}_I (M;\C)=1$, where $I$ is any of the three complex
structures appearing in $\goh$ and $H^{2,0}_I (M;\C)$ is the
$(2,0)$-part in the Hodge decomposition of $H^2 (M;\C)$ defined by
$I$ (for all three
complex structures in $\goh$ the space $H^{2,0}_\cdot (M;\C)$ has
the same dimension). K3-surfaces, as well as the Hilbert schemes of
points for K3-surfaces, are IHS. Any closed
hyperk\"ahler manifold admits a finite covering which is the product
of a torus and several IHS hyperk\"ahler manifolds
\cite{_Bogomolov:decompo_}. The IHS hyperk\"ahler manifolds are also
called {\bf hyperk\"ahler manifolds of maximal holonomy}, because
the holonomy group of a hyperk\"ahler manifold is $\Sp(n)$ (the
group of invertible quaternionic $n\times n$-matrices) -- and not
its proper subgroup -- if and only if it is IHS
\cite{_Besse:Einst_Manifo_}.

%%%%%%%%%%%%%%%%%%%%%%%%%%%%%%%%%%%%%%%%%%%%%%%%%%%%%%%%%%%%%%%%%%%%%%%

\hfill

%%%%%%%%%%%%%%%%%%%%%%%%%%%%%%%%%%%%%%%%%%%%%%%%%%%%%%%%%%%%%%%%%%%%%%%%

%%%%%%%%%%%%%%%%%%%%%%%%%%%%%%%%%%%%%%%%%%%%%%%%%%%%%%%%%%%%
\subsection{Symplectic packing of tori and IHS hyperk\"ahler manifolds by ellipsoids}
\label{_USP_tori_IHS_hyperk_mfds_Subsection_}
%%%%%%%%%%%%%%%%%%%%%%%%%%%%%%%%%%%%%%%%%%%%%%%%%%%%%%%%%%%%

By $\Vol$ we will always denote the symplectic volume of a symplectic manifold.

A {\bf closed ellipsoid} in $\C^n$ is defined as a
set \[\{ (z_1,\ldots,z_n)\in\C^n\ |\ \sum_{i=1}^n a_i |z_i|^2 \leq r\}\]
for some $a_1,\ldots, a_n, r >0$.

Let $(M,\omega)$, $\dim_\R M = 2n$, be a closed connected symplectic manifold.
We say that {\bf the symplectic packings of $(M,\omega)$ by ellipsoids are unobstructed},
if any finite collection of pairwise disjoint
closed ellipsoids in the standard symplectic $\R^{2n}$ of total volume less than
$\Vol (M,\omega)$ has an open neighborhood that can be symplectically embedded into $(M,\omega)$.

%%%%%%%%%%%%%%%%%%%%%%%%%%%%%%%%%%%%%%%%%%%%%%%%%%%%%%%%%%%%%%%%%%%%%%%%

\hfill

%%%%%%%%%%%%%%%%%%%%%%%%%%%%%%%%%%%%%%%%%%%%%%%%%%%%%%%%%%%%%%%%%%%%%%%%
\theorem\label{_USP_main_Theorem_}\\ Let $M$ be either a torus $T^{2n}$
with a K\"ahler form $\omega$ or an IHS hyperk\"ahler manifold with
a hyperk\"ahler symplectic form $\omega$. Then the symplectic packings of $(M, \omega)$ by ellipsoids
are unobstructed.

%%%%%%%%%%%%%%%%%%%%%%%%%%%%%%%%%%%%%%%%%%%%%%%%%%%%%%%%%%%%%%%%%%%%%%%%

\hfill

%%%%%%%%%%%%%%%%%%%%%%%%%%%%%%%%%%%%%%%%%%%%%%%%%%%%%%%%%%%%%%%%%%%%%%%%

\ref{_USP_main_Theorem_} follows from a similar result (see \ref{_USP_Campana_Simple_Theorem_}) for a wider class of K\"ahler manifolds as explained in Section~\ref{_Campana-simple-mfds_Section_} below.

%%%%%%%%%%%%%%%%%%%%%%%%%%%%%%%%%%%%%%%%%%%%%%%%%%%%%%%%%%%%%%%%%%%%%%%

\hfill

%%%%%%%%%%%%%%%%%%%%%%%%%%%%%%%%%%%%%%%%%%%%%%%%%%%%%%%%%%%%%%%%%%%%%%%%

%%%%%%%%%%%%%%%%%%%%%%%%%%%%%%%%%%%%%%%%%%%%%%%%%%%%%%%%%%%%
\subsection{Symplectic packing of arbitrary Campana simple manifolds}
\label{_Campana-simple-mfds_Section_}
%%%%%%%%%%%%%%%%%%%%%%%%%%%%%%%%%%%%%%%%%%%%%%%%%%%%%%%%%%%%

If $J$ is a complex structure of K\"ahler type on a closed connected
manifold $M$, then the union ${\goth U}$ of  all complex
subvarieties $Z\subset M$ satisfying $0< \dim_\C Z< \dim_\C M$
either has measure\footnote{The measure is defined by means of a volume form on $M$. One can easily see
that if a set is of measure zero with respect to some volume form, then it is of measure zero with respect
to any other volume form.} zero or is the whole $M$ (see \cite[Remark 4.2]{_EV-balls_}).

If ${\goth U}$ has measure zero, $J$ is called {\bf Campana simple} and
the points of $M\backslash {\goth U}$ are called {\bf  Campana-generic}.

We say that $(M,J)$ is a {\bf Campana simple} complex manifold, if $J$ is a Campana simple complex structure (of K\"ahler type) on $M$.

%%%%%%%%%%%%%%%%%%%%%%%%%%%%%%%%%%%%%%%%%%%%%%

\hfill

%%%%%%%%%%%%%%%%%%%%%%%%%%%%%%%%%%%%%%%%%%%%%%

\remark \\
Campana simple manifolds are non-algebraic.
According to a conjecture of Campana
(see \cite[Question 1.4]{_Campana:isotrivial_},
\cite[Conjecture 1.1]{_CDV:threefolds_}),
any Campana simple manifold
is bimeromorphic to a
hyperk\"ahler orbifold or a finite quotient of a torus.

%%%%%%%%%%%%%%%%%%%%%%%%%%%%%%%%%%%%%%%%%%%%%%

\hfill

%%%%%%%%%%%%%%%%%%%%%%%%%%%%%%%%%%%%%%%%%%%%%%

We say that a complex structure $J$ of K\"ahler type on $M$ {\bf can be approximated by Campana-simple complex
structures} (in a smooth deformation) if there exists a smooth family $\{ J_t\}_{t\in B^{2m}}$, $J_0=J$,
of complex structures $J_t$ on $M$
and a sequence $\{ t_i\}\to 0$ in $B^{2m}$ so that each $J_{t_i}$ is
Campana simple. (Here $B^{2m}\subset \C^m$ is an open ball centered at $0$). Note that it follows from a version of
Kodaira-Spencer stability theorem \cite{_Kod-Spen-AnnMath-1960_} (see \cite[Theorem 5.6]{_EV-balls_} for more details), that if $J$ is of K\"ahler type, then so are
$J_t$ for $t\in B^{2m}$ sufficiently close to $0$.

%%%%%%%%%%%%%%%%%%%%%%%%%%%%%%%%%%%%%%%%%%%%%%

\hfill

%%%%%%%%%%%%%%%%%%%%%%%%%%%%%%%%%%%%%%%%%%%%%%

\theorem\cite[Theorem 4.5]{_EV-balls_}\label{_Campana_simple_complex_structures_dense_for_tori_IHS_hyperk_mfds_Theorem_}\\
(A) Any complex structure of K\"ahler type on $T^{2n}$ can be approximated by
complex structures $J$
such that $(T^{2n},J)$ does not admit any proper complex subvarieties
of positive dimension and, in particular, is Campana simple.

\noindent (B) Let $(M,\goh)$, $\dim_\R M\geq 4$, be a closed connected IHS
hyperk\"ahler manifold and let $I$ be a complex structure appearing in $\goh$.
Then $I$ can be approximated by
Campana simple complex structures. \endproof

%%%%%%%%%%%%%%%%%%%%%%%%%%%%%%%%%%%%%%%%%%%%%%%

\hfill

%%%%%%%%%%%%%%%%%%%%%%%%%%%%%%%%%%%%%%%%%%%%%%%

\theorem\label{_USP_Campana_Simple_Theorem_}\\
Let $(M,I,\omega)$ be a K\"ahler manifold and assume that $I$ can be
approximated in a smooth deformation by Campana simple complex structures.
Then the symplectic packings of $(M,\omega)$ by ellipsoids are unobstructed.

%%%%%%%%%%%%%%%%%%%%%%%%%%%%%%%%%%%%%%%%%%%%%%%%%%%%%%%%%%%%%%%%%%%%%%%%

\hfill

%%%%%%%%%%%%%%%%%%%%%%%%%%%%%%%%%%%%%%%%%%%%%%%

For the proof of \ref{_USP_Campana_Simple_Theorem_} see Section~\ref{_USP_Campana_Simple_Theorem_Pf_Section_}.

%%%%%%%%%%%%%%%%%%%%%%%%%%%%%%%%%%%%%%%%%%%%%%%%%%%%%%%%%%%%%%%%%%%%%%%%

\hfill

%%%%%%%%%%%%%%%%%%%%%%%%%%%%%%%%%%%%%%%%%%%%%%%

\remark\\
In case $\dim_\R M = 4$ \ref{_USP_Campana_Simple_Theorem_} follows immediately from the analogous result for the symplectic packings by balls
proved in \cite{_EV-balls_} and an observation by McDuff \cite{_McDuff-ellipsoids_} that $(M,\omega)$ admits a symplectic embedding of a 4-dimensional ellipsoid if and only if
it admits a symplectic embedding of the disjoint union of a number of equal balls of the same total volume as the ellipsoid.
However, the proof of \ref{_USP_Campana_Simple_Theorem_} that we give below works in all dimensions.

%%%%%%%%%%%%%%%%%%%%%%%%%%%%%%%%%%%%%%%%%%%%%%%

\hfill

%%%%%%%%%%%%%%%%%%%%%%%%%%%%%%%%%%%%%%%%%%%%%%%

\noindent
{\bf Proof of \ref{_USP_main_Theorem_}.}\\ \ref{_USP_main_Theorem_} follows from \ref{_USP_Campana_Simple_Theorem_} and
\ref{_Campana_simple_complex_structures_dense_for_tori_IHS_hyperk_mfds_Theorem_}. \endproof

%%%%%%%%%%%%%%%%%%%%%%%%%%%%%%%%%%%%%%%%%%%%%%%
\subsection{Idea of the proof of \ref{_USP_Campana_Simple_Theorem_} and plan of the paper}
%%%%%%%%%%%%%%%%%%%%%%%%%%%%%%%%%%%%%%%%%%%%%%%

The idea of the proof of \ref{_USP_Campana_Simple_Theorem_} is as follows. For simplicity we outline it in the case when $(M,I)$ is Campana simple.

By a result on simultaneous Diophantine
approximation (see \ref{_approximation-by-pairwise-coprime-vectors_Proposition_}), any ellipsoid can be approximated by an ellipsoid
of the form $\{ \pi\sum_{i=1}^n a_i |z_i|^2 \leq r\}$ where all $a_i$ are pairwise coprime positive integers. (We call such an ellipsoid simple).
Thus, it suffices to prove that a disjoint union of $k$ simple ellipsoids of total volume less than $\Vol (M,\omega)$ admits a symplectic
embedding to $(M,\omega)$.

Recall that McDuff-Polterovich in \cite{_McD-Polt_}
have shown that the problem of symplectic packing
of a symplectic manifold $M$ by symplectic balls can be interpreted as a problem
about the shape of the symplectic cone of the symplectic blow-up of $M$.
There is a similar interpretation for symplectic packing by symplectic ellipsoids
where the role of the blow-ups is played by the weighted blow-ups (Subsection
\ref{_weighted-blowup_Section}).

By an extension of the McDuff-Polterovich results (for symplectic embeddings of balls) to the ellipsoid case,
it suffices to show that a certain degree-2 real cohomology class $\talpha$ of
the weighted blow-up $\tM$ of $M$ at some $k$ points $x_1,\ldots,x_k$ (with the weights given by the coefficients $a_i$ from the
equations of the ellipsoids) is K\"ahler. The weighted blow-up $\tM$ is a complex orbifold {\it with isolated singularities} due
to the fact that the ellipsoids are simple.
Moreover, by an extension of a result of McDuff-Polterovich (for symplectic embeddings of balls) to the ellipsoid case, the
complex orbifold $\tM$ admits
a K\"ahler structure.
We construct a K\"ahler resolution $\pi:\hN\to \tM$ (our construction uses the fact that $\tM$ has only isolated singularities)
and consider a cohomology class of the form $\pi^* \talpha + \delta b$, where $\delta>0$ is small and $b\in H^2 (\hN;\R)$
has the property $\pi_* b =0$ and $b\cup \pi^* \talpha =0$. Using the Demailly-Paun theorem, describing the K\"ahler cone of $\hN$,
and the fact that the points $x_1,\ldots,x_k$ are Campana generic we show that for a sufficiently small $\delta>0$ the class
$\pi^* \talpha + \delta b$ is K\"ahler. Then, using a result of Miyaoka on the extension of a K\"ahler form over an isolated puncture,
we show that the class $\pi_* (\pi^* \talpha + \delta b)= \talpha$ is K\"ahler which finishes the proof.

The plan of the paper is as follows.

In Section~\ref{_orbifolds-definition_Section_} we recall basic facts about orbifolds.

In Section~\ref{_resolutions_of_orbifolds_Section_} we state the result on the existence of a
K\"ahler resolution of a closed K\"ahler orbifold with isolated singularities and prove that the pushforward of a K\"ahler class
of the resolution is a K\"ahler class of the base.

In Section~\ref{_construction_Kahler_resolution_Section_} we construct the resolution and prove the existence of a K\"ahler form on it. The construction uses
plurisubharmonic functions -- we recall basic facts about them in the beginning of the same section.

In Section~\ref{_weighted-blowup_Section} we recall the basics concerning weighted blow-ups in the complex and symplectic category
and discuss the relation between the symplectic/K\"ahler classes on weighted blowups and symplectic embeddings of ellipsoids.

In Section~\ref{_Demailly_Paun_Section_} we apply the Demailly-Paun theorem to prove that the needed cohomology class of the resolution
is K\"ahler.

In Section~\ref{_USP_Campana_Simple_Theorem_Pf_Section_} we combine the previous results together and prove \ref{_USP_Campana_Simple_Theorem_}.

In the appendix (Section~\ref{_approximation_appendix_Section_}) we prove the above-mentioned result on the simultaneous Diophantine approximation.

%%%%%%%%%%%%%%%%%%%%%%%%%%%%%%%%%%%%%%%%%%%%%%%%%%%%%%%%%%%%%%%%%%%%%%%%

\section[Orbifolds, weighted blow-ups and symplectic packing  by ellipsoids]%
{Orbifolds, weighted blow-ups \\ and symplectic packing by ellipsoids}
\label{_orbifolds_Section_}
%%%%%%%%%%%%%%%%%%%%%%%%%%%%%%%%%%%%%%%%%%%%%%%%%%%%%%%%%%%%%%%%%%%%%%%%

Here we recall a few basic facts on orbifolds.
Orbifolds were originally introduced under the name ``V-manifolds" by I.Satake in 1950s \cite{_Satake-PNAS-1956_}.
The name was changed to ``orbifolds" in W.Thurston's seminar in the 1970s \cite{_Thurston-MSRI2002_}.

%%%%%%%%%%%%%%%%%%%%%%%%%%%%%%%%%%%%%%%%%%%%%%%%%%%%%%%%%%%%%%%%%%%%%%%%

\subsection{Basics of orbifolds}
\label{_orbifolds-definition_Section_}
%%%%%%%%%%%%%%%%%%%%%%%%%%%%%%%%%%%%%%%%%%%%%%%%%%%%%%%%%%%%%%%%%%%%%%%%

\hfill

\definition
A {\bf real (complex) orbifold chart} (also known as a {\bf locally uniformizing system}) on a topological space $N$ consists of the following objects:

\smallskip
\noindent
$U_\alpha$ -- an open connected subset of $N$.

\smallskip
\noindent
$\tU_\alpha$ -- an open connected neighborhood of the origin $0$ in $\R^n$ (in $\C^n$).

\smallskip
\noindent
$\Gamma_\alpha$ -- a finite (possibly trivial) group acting effectively on $\R^n$ (on $\C^n$) by linear
real (complex) transformations so that $\tU_\alpha$ is invariant under the action,
$0$ is a fixed point of the action and the set of points of $\tU_\alpha$ with a non-trivial stabilizer
is of real codimension 2 (complex codimension 1) or more.

\smallskip
\noindent
$\phi_\alpha: U_\alpha\to\tU_\alpha/\Gamma_\alpha$ -- a homeomorphism.

\smallskip
The number $n$ is called the real (complex) {\bf dimension of the chart}. If $x\in U_\alpha$, we call the stabilizer of a pre-image of $\phi_\alpha (x)\in \tU_\alpha/\Gamma_\alpha$ in $\tU_\alpha$ under the action of $\Gamma_\alpha$ {\bf the stabilizer of $x$ in $U_\alpha$} and denote it by $\Gamma_{\alpha,x}$.

\bigskip
A {\bf real (complex) orbifold atlas} on $N$ is a collection
\[ \{ (U_\alpha,\tU_\alpha, \Gamma_\alpha, \phi_\alpha)\}\]
of $n$-dimensional
real (complex) orbifold charts on $N$ with the following properties:

\begin{itemize}

\item{} $N=\cup_\alpha U_\alpha$.

\item{} Any finite intersection of sets from the collection $\{ U_\alpha\}$ is a union of sets from the collection.

\item{} If $U_\alpha\subset U_\beta$, then there exists
an injective homomorphism $f_{\alpha\beta}:\Gamma_\alpha\to\Gamma_\beta$
a smooth (complex analytic) embedding $\tphi_{\alpha\beta}:\tU_\alpha\to \tU_\beta$ equivariant with respect to
the actions of $\Gamma_\alpha,\Gamma_\beta$ (related by $f_{\alpha\beta}$) and covering the inclusion
$U_\alpha\hookrightarrow U_\beta$ (where $U_\alpha, U_\beta$ are identified with $\tU_\alpha/\Gamma_\alpha$ and $\tU_\beta/\Gamma_\beta$, respectively, by $\phi_\alpha$ and $\phi_\beta$).

\end{itemize}

\bigskip
A {\bf real (complex) orbifold of real (complex) dimension $n\geq 2$} is a Hausdorff paracompact topological space $N$
equipped with a maximal real (complex) orbifold atlas formed by $n$-dimensional orbifold charts. Such a maximal atlas is called an {\bf orbifold structure}.

\bigskip
Given a point $x$ in an orbifold $N$, the stabilizers $\Gamma_{\alpha,x}$ of $x$ in different orbifold charts $U_\alpha$ on $N$ containing $x$ are all isomorphic. We will denote any of these stabilizers by $\Gamma_x$, or by $\Gamma_x^N$, if we want to emphasize that it is the stabilizer of $x$ in $N$.

The {\bf tangent space}
of an $n$-dimensional orbifold $N$ at $x\in N$ is defined as the $n$-dimensional {\sl representation} of $\Gamma_x \cong \Gamma_{\alpha,x}$ on $T_{\phi_\alpha (x)} \tU_\alpha$, where
$U_\alpha$ is an orbifold chart containing $x$. The representations of $\Gamma_x$ coming from different orbifold charts $U_\alpha, U_\beta$ containing $x$ are isomorphic (as representations of isomorphic groups $\Gamma_{\alpha,x}\cong\Gamma_{\beta,x}$).

If the stabilizer of a point in $N$ is trivial, then it is called a {\bf regular point} of $N$; otherwise it is
called a {\bf singular point}.
The set of singular points of $N$
is called {\bf the singular locus of $N$} and its complement {\bf the regular part} of $N$.
The regular part of $N$ is a smooth manifold -- it admits a smooth atlas
formed by some of the orbifold charts from the maximal orbifold atlas on $N$.
We say that $N$ is {\bf an orbifold  with isolated singularities} if the singular locus of $N$ is a discrete set.

Of course, any smooth real (complex) manifold is also a real (complex) orbifold. We will say that a real (complex) orbifold is {\bf smooth}, if its orbifold structure contains a smooth real (complex) subatlas or, equivalently, its singular locus is empty.

Most differential-geometric objects (smooth/complex analytic maps and their differentials, vector fields and their flows, differential forms, the differential of a differential form, Lie derivative along a vector field, almost complex structures, Riemannian metrics, vector bundles etc.)
can be generalized to orbifolds in a straightforward way: first, one considers a $\Gamma_\alpha$-invariant (or equivariant) version of an object on each $\tU_\alpha$ and then requires that the maps $\tphi_{\alpha\beta}$ glue the objects on all $\tU_\alpha$ and $\tU_\beta$.
In order to distinguish between the objects on smooth manifolds and their counterparts on orbifolds we will use the prefix ``orbifold" for the latter counterparts: orbifold smooth functions, orbifold smooth vector fields etc. Note that the restriction of an orbifold smooth function, orbifold smooth vector field etc. on an orbifold $N$ to the regular part $N^{reg}$ of $N$ is a smooth function, smooth vector field etc. in the usual sense on the smooth manifold $N^{reg}$.

Orbifolds admit partitions of unity (see e.g \cite[Theorem B.12]{_Barletta-Dragomir-Duggal_}). This allows to equip any orbifold with an orbifold Riemannian metric and to define the integral of an orbifold differential form over an oriented orbifold (use a partition of unity subordinated to orbifold charts and for an $n$-dimensional chart $\{ (U_\alpha,\tU_\alpha, \Gamma_\alpha, \phi_\alpha)\}$ and an orbifold form $\Omega$ of degree $n$ supported in $U_\alpha$ define $\int_{U_\alpha} \Omega$ as $1/|\Gamma_\alpha| \int_{\tU_\alpha} \tOmega$, where $\tOmega$ is the lift of $\Omega$ to $\tU_\alpha$).

A suborbifold $L$ of a real (complex) orbifold $N$ is then defined
as a subset $L\subset N$ equipped with an orbifold structure so that
the inclusion $L\hookrightarrow N$ is an orbifold smooth (analytic) map.
In particular, this means that $\Gamma_x^L$
injects into $\Gamma_x^N$ for every $x\in L$. A smooth suborbifold
of an orbifold $N$ will be called a {\bf submanifold} of $N$.

%%%%%%%%%%%%%%%%%%%%%%%%%%%%%%%%%

\hfill

%%%%%%%%%%%%%%%%%%%%%%%%%%%%%%%%

\definition

A {\bf pairwise coprime vector} is an ordered tuple of pairwise coprime positive integers.

Given a vector $\ba=(a_1,\ldots,a_n)$ set
\[
\langle\bx\rangle:=a_1\cdot\ldots\cdot a_n.
\]

%%%%%%%%%%%%%%%%%%%%%%%%%%%%%%%%%

\hfill

%%%%%%%%%%%%%%%%%%%%%%%%%%%%%%%%

\example \label{_weighted-proj-space_Example_}

Let $\ba:=(a_0, \ldots,a_n)$ be a pairwise coprime vector. The {\bf weighted
projective space} $\C P^n (a_0,\ldots, a_n)$ (which we will also
denote by $\C P^n(\ba)$) is defined as the quotient of $\C^{n+1}$ by
the action of $\C^*$ given by
\begin{equation}
\label{_eqn-Cstar-action_} \lambda:
(z_0, \ldots, z_n) \mapsto (\lambda^{a_0} z_0, \ldots,
\lambda^{a_n}z_n).
\end{equation}

The weighted projective space $\C P^n (\ba)$ can be equipped with
the structure of a complex orbifold (see e.g. \cite{_Godinho_}). Since
$a_0,\ldots,a_n$ are pairwise coprime, the singular
locus of $\C P^n (\ba)$ is discrete.

The integral homology/cohomology of $\C P^n (\ba)$ is isomorphic (as
a group) to that of $\C P^n$, while the
multiplicative structure of $H^* (\C P^n (\ba);\Q)$ may differ from
the case of $\C P^n$ (see \cite{_Tramer_},\cite{_Kawasaki_}) and does depend
on $\ba$. In particular, in our case, when all $a_0,\ldots,a_n$ are
pairwise coprime, the map
\[
\C P^n\to \C P^n (\ba),
\]
given by
\[
[z_0:\ldots:z_n]\mapsto [z_0^{a_0}:\ldots:z_n^{a_n}],
\]
induces a ring isomorphism
\[
H^* (\C P^n (\ba);\Z)\to H^* (\C P^n ;\Z)
\]
which is the multiplication by $\langle\ba\rangle$ in each degree. For each $i=0,\ldots,n$ denote by $\alpha_i$ the generator of $H^{2i} (\C P^n
(\ba);\Z)$ mapped by this ring isomorphism to the positive multiple of the
standard generator of $H^{2i} (\C P^n;\Z)$. Then $\alpha_i \alpha_j = \langle\ba\rangle\alpha_{i+j}$ for all $i,j$, $i+j\leq n$. In particular, $\alpha_1^n = \langle\ba\rangle^{n-1}\alpha_n$.

The hyperplane section $z_0=0$ is a suborbifold $L$ of $\C P^n (\ba)$ which is orbifold
diffeomorphic to $\C P^{n-1} (a_1,\ldots, a_n)$. Then $\alpha_1\in H^2 (\C P^n (\ba);\Z)$ is
Poincar\'e-dual to $[L]$ and
\[
\langle \alpha_1^n, [\C P^n (\ba)]\rangle = \langle\ba\rangle^{n-1}.
\]

The space $\C P^n (\ba)$ can be equipped with an orbifold symplectic
structure. Namely, consider the Hamiltonian
$H(z_0,\ldots,z_n) = \pi \sum_{i=0}^n a_i |z_i|^2$ on the standard
symplectic $\C^{n+1}$. It defines an $S^1=\R/\Z$-action:
\begin{equation}
\label{_eqn-S1-action_} t: (z_0, \ldots, z_n) \mapsto (e^{2\pi\1 a_0
t}z_0, \ldots, e^{2\pi\1 a_n t}z_n).
\end{equation}
which is the restriction of the $\C^*$-action
\eqref{_eqn-Cstar-action_} to $S^1\subset \C^*$. For $r>0$ consider
the reduced space $H^{-1} (r)/S^1$ -- it has the structure of a real
$2n$-dimensional orbifold (see e.g. \cite{_Godinho_}) and is
naturally identified (by an orbifold diffeomorphism) with $\C
P^n(\ba)$. The reduction induces an orbifold symplectic form on $H^{-1}
(r)/S^1$ and hence on $\C P^n (\ba)$. We will denote this orbifold form by
$\Omega_{\ba,r}$. Set $\Omega_\ba := \Omega_{\ba,1}$ -- we will call it {\bf the Fubini-Study symplectic form
on $\C P^n (\ba)$}. One can check that the form $\Omega_\ba$ is K\"ahler, $\Omega_{\ba,r} = r\Omega_\ba$ and
\[
[\Omega_{\ba,r}]=
\frac{r}{\langle\ba\rangle} \alpha_1,
\]
\[
\int_{\C P^n (\ba)} \Omega_{\ba,r}^n =
\frac{r^n}{\langle\ba\rangle} \Vol_{2n},
\]
where $\Vol_{2n}$ is the volume of the Euclidean $2n$-dimensional unit
ball.

%%%%%%%%%%%%%%%%%%%%%%%%%%%%%%%%%%%%%%%%%%%%%%%%%%%%%%%%%%%%

\hfill

Poincar\'e duality (over $\Q$),
de Rham and Hodge theorems for closed manifolds extend to closed (=compact, without boundary) orbifolds --
see \cite{_Satake-PNAS-1956_} and \cite{_Baily-AJM-1956_}.
This allows to obtain an orbifold version of Moser's theorem \cite{_Moser_} -- its proof literally repeats the proof for the smooth case and we
write it here as an example of a straightforward generalization of a result for smooth manifolds to orbifolds.

\hfill

%%%%%%%%%%%%%%%%%%%%%%%%%%%%%%%%%%%%%%%%%%%%%%%%%%%%%%%%%%%%

\theorem
\label{_Moser_Proposition_}\\
Let $N$ be a closed symplectic orbifold, and
$\omega_t$ a smooth family of orbifold symplectic forms, parameterized by $t\in [0,1]$
and lying in the same (de Rham) cohomology class.
Then there exists a smooth family of orbifold diffeomorphisms $\Psi_t:\ N \arrow N$
such that $\Psi_t^*\omega_0 = \omega_t$.

%%%%%%%%%%%%%%%%%%%%%%%%%%%%%%%%%%%%%%%%%%%%%%%%%%%%%%%%%%%%

\hfill

\noindent
{\bf Proof.}\\
The orbifold form $\dot \omega_t$ is exact
and therefore there exists an orbifold
1-form $\eta_t$ satisfying
$d\eta_t= \dot \omega_t$. The Hodge theorem allows to choose the orbifold forms $\eta_t$ so that they depend on $t$ smoothly.
Let $X_t$ be an orbifold vector field satisfying $\omega_t \cntrct X_t=\eta_t$.
Then $\Lie_{X_t}\omega_t = \eta_t=\dot \omega_t$ by Cartan's formula (which also extends to orbifolds).
The orbifold vector field $X_t$ defines a flow of orbifold diffeomorphisms $\Psi_t$ on $N$.
Then $\Lie_{X_t}\omega_t = \dot \omega_t$ implies
$\Psi_t^*(\omega_0)=\omega_t$.
\endproof

%%%%%%%%%%%%%%%%%%%%%%%%%%%%%%%%%%%%%%%%%%%%%%%%%%%%%%%%%%%%

\hfill

\remark
\label{_Moser_local_versions_Remark_}\\
Similarly one can prove orbifold versions of various symplectic neighborhood theorems -- they are all based on Moser's method as above
and easily generalize to orbifolds.

\hfill

%%%%%%%%%%%%%%%%%%%%%%%%%%%%%%%%%%%%%%%%%%%%%%%%%%%%%%%%%%%%%%%%%%%%%%%%

\subsection{Resolution of orbifolds and K\"ahler classes}
\label{_resolutions_of_orbifolds_Section_}
%%%%%%%%%%%%%%%%%%%%%%%%%%%%%%%%%%%%%%%%%%%%%%%%%%%%%%%%%%%%%%%%%%%%%%%%

Given a smooth map $F: M\to N$ between closed oriented smooth manifolds, we define the pushforward $F_* b\in H^* (N;\R)$
of a cohomology class $b\in H^* (M;\R)$ using the Poincar\'e-duality on $M$ and $N$ and the pushforward of homology classes.

Similarly to the smooth case, a real $(1,1)$-cohomology class of a closed complex orbifold $(N,J)$ is called {\bf K\"ahler} if it can be realized by an orbifold K\"ahler form. Denote by $\Kah (N,J)$ the cone in $H^2 (N;\R)$ formed by the K\"ahler cohomology classes.

%%%%%%%%%%%%%%%%%%%%%%%%%%%%%%%%%%%%%%%%%%%%%%%%%%%%%%%%%%%%%%%%%%%%%%%%

\hfill

%%%%%%%%%%%%%%%%%%%%%%%%%%%%%%%%%%%%%%%%%%%%%%%%%%%%%%%%%%%%%%%%%%%%%%%%

\theorem\label{_pushforward_of_Kahler_classes_resolution_of_orbifolds_Theorem_}\\
Let $N$, $\dim_\C N = n\geq 2$, be a closed complex orbifold with isolated singularities. Denote its singular locus by $\Sigma = \{ y_1,\ldots,y_m\}$.

\smallskip
\noindent
A. Let $N'$, $\dim_\C N = n$, be a closed complex manifold and
let $P: N'\to N$ be a surjective smooth map such that
$P: N'\setminus P^{-1} (\Sigma)\to N\setminus \Sigma$ is a biholomorphism and $P^{-1} (\Sigma)$ has zero volume
(with respect to a volume form on $N'$).

Then the pushforward of any K\"ahler cohomology class on $N'$
is a K\"ahler cohomology class on $N$.

\smallskip
\noindent
B. There exist

\smallskip
\noindent
- a smooth closed complex manifold $\hN$ of the same dimension as $N$;

\noindent
- a holomorphic map $\pi: \hN\to N$ such that $\pi: \hN\setminus \pi^{-1} (\Sigma)\to N\setminus\Sigma$ is a biholomorphism;

\noindent
- cohomology classes $b_i\in H^2 (\hN;\R)$, $i=1,\ldots,m$, so that $\pi_* b_i = 0$ and $b_i \cup \pi^* u =0$ for any $i$ and any $u\in H^* (N)$, $\textrm{deg}\ u >0$;

\smallskip
\noindent
so that for any $v\in \Kah (N,J)$ and any sufficiently small $\delta_1,\ldots,\delta_m>0$,
\[
\pi^* v + \sum_{i=1}^m \delta_i b_i \in \Kah (\hN,\hJ).
\]

%%%%%%%%%%%%%%%%%%%%%%%%%%%%%%%%%%%%%%%%%%%%%%%%%%%%%%%%%%%%%%%%%%%%%%%%

\hfill

%%%%%%%%%%%%%%%%%%%%%%%%%%%%%%%%%%%%%%%%%%%%%%%%%%%%%%%%%%%%%%%%%%%%%%%%

The construction of $\hN$ and $\pi$ in part B of \ref{_pushforward_of_Kahler_classes_resolution_of_orbifolds_Theorem_} amounts to resolving
the isolated singularities of a K\"ahler orbifold $N$ in the K\"ahler category.
In particular, this yields the following corollary.

%%%%%%%%%%%%%%%%%%%%%%%%%%%%%%%%%%%%%%%%%%%%%%%%%%%%%%%%%%%%%%%%%%%%%%%%

\hfill

%%%%%%%%%%%%%%%%%%%%%%%%%%%%%%%%%%%%%%%%%%%%%%%%%%%%%%%%%%%%%%%%%%%%%%%%

\corollary Any closed K\"ahler orbifold with isolated singularities admits a K\"ahler resolution. \endproof

%%%%%%%%%%%%%%%%%%%%%%%%%%%%%%%%%%%%%%%%%%%%%%%%%%%%%%%%%%%%%%%%%%%%%%%%

\hfill

%%%%%%%%%%%%%%%%%%%%%%%%%%%%%%%%%%%%%%%%%%%%%%%%%%%%%%%%%%%%%%%%%%%%%%%%

Part A of \ref{_pushforward_of_Kahler_classes_resolution_of_orbifolds_Theorem_} is proved below.

For the proof of part B of \ref{_pushforward_of_Kahler_classes_resolution_of_orbifolds_Theorem_} see Section~\ref{_construction_Kahler_resolution_Section_}.

%%%%%%%%%%%%%%%%%%%%%%%%%%%%%%%%%%%%%%%%%%%%%%%%%%%%%%%%%%%%%%%%%%%%%%%%

\hfill

%%%%%%%%%%%%%%%%%%%%%%%%%%%%%%%%%%%%%%%%%%%%%%%%%%%%%%%%%%%%%%%%%%%%%%%%

\noindent
{\bf Proof of
part A of \eqref{_pushforward_of_Kahler_classes_resolution_of_orbifolds_Theorem_}}.

Let $\theta$ be a K\"ahler form on $N'$. Then $P_* \theta$ is a K\"ahler form on the smooth complex manifold $N\setminus \Sigma$.

Assume $x\in \Sigma$ and $\Gamma_x$ is its stabilizer. A neighborhood $U$ of $x$ in $N$ is biholomorphic to a neighborhood of zero in $\C^n/\Gamma_x$ and, since $x$ is an isolated singularity of $N$, the form $P_* \theta|_{U\setminus x}$ lifts under the projection $\C^n\to \C^n/\Gamma_x$ to a $\Gamma_x$-invariant K\"ahler form $\zeta$ on $V\setminus 0$, where $V\subset \C^n$ is a $\Gamma_x$-invariant open set which is the lift of $U$.
By a result of \cite{_Miyaoka_}, there exists a
K\"ahler form $\zeta'$ on the whole $V$ that coincides with $\zeta$ outside a small neighborhood of $0$. Averaging, if necessary, $\zeta'$ with respect to the
action of $\Gamma_x$ on $V$, we can assume that $\zeta'$ is a $\Gamma_x$-invariant. Thus, $\zeta'$ descends to an orbifold K\"ahler form on
 $U$ that coincides with $P_* \theta$ outside a small neighborhood of $x$ in $U$. Thus, $P_* \theta$ can be extended from
 $N\setminus \Sigma$ to an orbifold K\"ahler form on $N$. This K\"ahler form coincides with $P_* \theta$ outside a finite union of disjoint contractible sets and therefore its cohomology class is $P_* [\theta]$. \endproof

\hfill

%%%%%%%%%%%%%%%%%%%%%%%%%%%%%%%%%%%%%%%%%%%%%%%%%%%%%%%%%%%%%%%%%%%%%%%%

%%%%%%%%%%%%%%%%%%%%%%%%%%%%%%%%%%%%%%%%%%%%%%%%%%%%%%%%%%%%%%%%%%%%%%%%

\subsection{Construction of a K\"ahler resolution}
\label{_construction_Kahler_resolution_Section_}
%%%%%%%%%%%%%%%%%%%%%%%%%%%%%%%%%%%%%%%%%%%%%%%%%%%%%%%%%%%%%%%%%%%%%%%%

We recall a few basic facts about currents and plurisubharmonic functions needed for the proof of \ref{_pushforward_of_Kahler_classes_resolution_of_orbifolds_Theorem_}. For more details see e.g. \cite{_Demailly_Compl_Analytic_Diff_Geom_,_Grif-Har_, _Lelong-Gruman_}.

Recall that a function $\varphi$, with values in $\R\cup \{ -\infty\}$, on an open domain $U\subset \C^n$ is called plurisubharmonic if it is upper semi-continuous (hence, locally bounded from above), not identically equal to $-\infty$ on any open set, and for any complex line $L$ in $\C^n$ the restriction of $\varphi$ to $U\cap L$ is either subharmonic or identically equal to $-\infty$. Plurisubharmonic functions are locally integrable \cite[Proposition I.9]{_Lelong-Gruman_}. A function $\varphi: U\to \R\cup \{ -\infty\}$ is called strictly plurisubharmonic if for every $p\in U$ and any sufficiently small
$\epsilon>0$ the function $\varphi-\epsilon |z|^2$ is plurisubharmonic. (Here and below $z=(z_1,\ldots, z_n)$ and $|z|^2:= |z_1|^2+\ldots+|z_n|^2$).
A smooth function $\varphi$ is (strictly) plurisubharmonic if and only if the $(1,1)$-form $\sqrt{-1}\partial\bar{\partial} \varphi$ is a (strictly) positive Hermitian form (being strictly positive is equivalent for the form $\sqrt{-1}\partial\bar{\partial} \varphi$ to being K\"ahler).

If $\varphi_1, \varphi_2$ are two -- not necessarily smooth -- (strictly) plurisubharmonic functions on $U\subset \C^n$,
then for any $\epsilon>0$ there exists
(see \cite[Lemma 5.18]{_Demailly_Compl_Analytic_Diff_Geom_})
a (strictly) plurisubharmonic function $\max_\epsilon \{ \varphi_1, \varphi_2\}$ on $U$, called the regularized maximum of $\varphi_1,\varphi_2$, such that

\begin{itemize}

\item{} if $\varphi_1, \varphi_2$ are smooth near $x\in U$, then so is $\max_\epsilon \{ \varphi_1, \varphi_2\}$;

\item{}
$\max \{ \varphi_1,\varphi_2\}
\leq
\max_\epsilon \{ \varphi_1, \varphi_2\}
\leq
\max\{ \varphi_1, \varphi_2\} + \epsilon$ on $U$;

\item{}
$\max_\epsilon \{ \varphi_1, \varphi_2\} (x) = \varphi_j (x)$, if $\varphi_i (x) +2\epsilon \leq \varphi_j (x)$, $i,j=1,2$, $i\neq j$, $x\in U$.

\end{itemize}

A current of degree $(p,q)$ on a complex manifold of complex
dimension $n$ is a continuous linear functional on the space of smooth
compactly supported differential complex-valued $(2n-p-q)$-forms that
vanishes on $(k,l)$-forms as long as $(k,l)\neq (n-p,n-q)$.

A current $T$ is called real if $\overline{T(\xi)} = T(\overline{\xi})$ for all differential forms $\xi$ (the bar is the complex conjugation).

Each differential $(p,q)$-form $\theta$ with locally integrable coefficients on a complex manifold defines a current $T_\theta$ of degree $(p,q)$: the value of $T_\theta$ on a smooth compactly supported $(2n-p-q)$-form $\xi$ is defined as the integral of $\theta\wedge\xi$ over the manifold. If two differential forms with continuous coefficients define the same current they coincide everywhere. If $\theta$ is real, then so is $T_\theta$. The currents defined in this way by smooth forms $\theta$ are called smooth.

A real current $T$ of degree $(p,p)$ is called positive if $T(\eta\wedge \bar{\eta})\geq 0$ for any
smooth compactly supported differential $(n-p,0)$-form, and strictly positive if $T-T'$ is positive for some positive current $T'$ of degree $(p,p)$.

Recall that a real differential $(1,1)$-form $\theta$ is called positive (strictly positive) if $\theta (v,\bar{v}) \geq 0$ ($\theta(v,\bar{v}) >0$) for any tangent vector $v\neq 0$. Thus, a positive (strictly positive) real differential form $(1,1)$-form $\theta$ defines a positive (strictly positive) real current $T$ of degree $(1,1)$.

The differentials $d, \partial, \bar{\partial}$ for currents on a complex manifold are defined by duality using the corresponding differentials for smooth complex-valued differential forms. The homomorphism $\theta\mapsto T_\theta$ induces a homomorphism between the cochain complexes of differential forms and currents (for any of the differentials $d, \partial, \bar{\partial}$) that induces an isomorphism between the corresponding cohomologies (see e.g. \cite[p.385]{_Grif-Har_}).

A smooth real $(1,1)$-form $\theta$ is K\"ahler if and only if the smooth current $T_\theta$ is closed and strictly positive.

By the local $\partial\bar{\partial}$-lemma \cite[Theorem 2.28]{_Lelong-Gruman_}, any (not necessarily smooth) closed real (strictly) positive current $T$ of degree $(1,1)$ on a complex manifold can be represented locally as $T=\sqrt{-1}\partial\bar{\partial} T_\varphi$ for a (strictly) plurisubharmonic function $\varphi$. The function $\varphi$ is defined uniquely up to an addition of a harmonic (hence, smooth) function. Thus, if $T$ is smooth near a point then $\varphi$ is also smooth near that point.

%%%%%%%%%%%%%%%%%%%%%%%%%%%%%%%%%%%%%%%%%%%%%%%%%%%%%%%%%%%%%%%%%%%%%%%%

\hfill

%%%%%%%%%%%%%%%%%%%%%%%%%%%%%%%%%%%%%%%%%%%%%%%%%%%%%%%%%%%%%%%%%%%%%%%%

\noindent
{\bf Proof of
part B of \ref{_pushforward_of_Kahler_classes_resolution_of_orbifolds_Theorem_}}.

Let $\Gamma_1,\ldots,\Gamma_m$ be the stabilizers of the singular points $y_1,\ldots,y_m$ of $N$.
Consider orbifold charts $(U_i,\tU_i, \Gamma_i, \phi_i)$ on $N$ such that $y_i\in U_i$ for all $i$ and all $U_i$ are pairwise disjoint.

Pick an arbitrary $i=1,\ldots,m$.
The singular point $y_i$ is isolated and therefore for a smaller neighborhood $U'_i\subset U_i$ of $y_i$ the punctured neighborhood $U'_i\setminus y_i$ can be bi-holomorphically identified with  $B(R)\setminus 0$, where $B (R)\subset \C^n$ denotes an open ball of radius $R$ centered at zero.

The orbifold $\C^n/\Gamma_i$ is a quasi-projective
algebraic variety with a single singularity at the
origin. By the famous result of Hironaka
\cite{_Hironaka_} (see also \cite{_Bierstone-Milman_}, \cite{_Villamayor1_,_Villamayor2_} for more accessible proofs of Hironaka's result), there exists a resolution $\pi_i:
X_i\to \C^n/\Gamma_i$ of the singularity where $X_i$ is a
smooth quasi-projective variety. Moreover, the resolution of
singularity $\pi_i: X_i\to \C^n/\Gamma_i$ is biholomorphic outside
of the singular set of $X_i$.
Define $V_i: = \pi_i^{-1} (\phi_i (U'_i))\subset X_i$. Since $\pi_i: X_i\setminus \pi_i^{-1} (0)\to \C^n/\Gamma_i\setminus 0$ is a biholomorphism,  the maps $\pi_i^{-1}$ and, respectively, $\pi_i^{-1}\circ \phi_i$ identify $B(R)\setminus 0\subset \C^n$ and, respectively, $U'_i\setminus y_i\subset N$ biholomorphically with the same open subset $W_i:= V_i\setminus \pi_i^{-1} (0)$.

Let us attach all the $V_i$ to $N\setminus \Sigma$ using the identifications $\pi_i^{-1}\circ \phi_i$ between the open sets
$U'_i\setminus y_i\subset N$ and $W_i= V_i\setminus \pi_i^{-1} (0)$.
As a result we get a complex manifold $\hN:= (N\setminus \Sigma)\cup \bigcup_{i=1}^m V_i$ along with a holomorphic projection $\pi: \hN\to N$ which is a biholomorphism over $\hN\setminus \pi^{-1} (\Sigma)$ and which coincides with $\pi_i$ over each $V_i$.

Assume $\eta$ is an orbifold K\"ahler form on $N$ such that $[\eta]=v$. The restriction of $\eta$ to the smooth complex manifold $N\setminus\Sigma$ can be viewed as a usual smooth K\"ahler
form on a smooth complex manifold.
In particular, under the identification $\phi_i$ the form $\eta|_{U'_i\setminus y_i}$ is identified with a smooth K\"ahler form on
$B(R)\setminus 0$ that will be also denoted by $\eta$.

As a smooth quasi-projective variety, $X_i$ carries a K\"ahler form $\omega_i$
(induced by the Fubini-Study form on the projective space of which $X_i$ is a subvariety).
By means of the identifications above, $\omega_i$ induces K\"ahler forms on $B(R)\setminus 0$ and $U'_i\setminus y_i$ that, by an abuse of notation, will
be both denoted also by $\omega_i$.

%%%%%%%%%%%%%%%%%%%%%%%%%%%%%%%%%%%%%%%%%%%%%%%%%

\hfill

%%%%%%%%%%%%%%%%%%%%%%%%%%%%%%%%%%%%%%%%%%%%%%%%%

\lemma
\label{lem-gluing-two-psh-functions}\\
For any sufficiently small $\delta_i >0$ there exists a K\"ahler form
$\xi_{i,\delta_i}$
on $B(R)\setminus 0$
that equals to $\delta_i \omega_i$ near $0$ and to $\eta$ outside $B(3R/4)$.

%%%%%%%%%%%%%%%%%%%%%%%%%%%%%%%%%%%%%%%%%%%%%%%%%

\hfill

%%%%%%%%%%%%%%%%%%%%%%%%%%%%%%%%%%%%%%%%%%%%%%%%%

Postponing the proof of the lemma let us finish the proof of the theorem.

By means of the identification $\phi_i$, the form $\xi_{i,\delta_i}$ induces a K\"ahler form on $U'_i\setminus y_i$, that will be also denoted by $\xi_{i,\delta_i}$, which coincides with
$\eta$ near the boundary of $U'_i$ and with $\omega_i$ near $y_i$.

It follows that for any sufficiently small $\delta_1,\ldots,\delta_m >0$, the manifold $\hN$ carries a K\"ahler form $\heta$ which, by definition, is equal to $\eta$ on $N\setminus \cup_{i=1}^m U'_i$ (identified by $\pi$ with a subset of $\hN$), to $\xi_{i,\delta_i}$ on each $U'_i$ (identified by $\pi_i^{-1}\circ \phi_i$ with $W_i$) and to $\delta_i\omega_i$ on each $V_i$.

The closed 2-forms $\heta$ and $\pi^* \eta$ coincide outside $\cup_{i=1}^m V_i$. Hence, the form $\heta - \pi^* \eta$ is a sum of closed 2-forms $\delta_i\sigma_i$, $i=1,\ldots,m$, on $\hN$ so that each $\sigma_i$ is K\"ahler, supported inside $V_i$ (and therefore $\pi_* [\sigma_i]=0$) and coincides with the form $\omega_i$ on the analytic subvariety $\pi^{-1} (y_i)\subset V_i$. Each form $\omega_i$ is K\"ahler on $V_i$ and $\pi^{-1} (y_i)$ is a deformation retract of $V_i$ -- therefore the cohomology class $b_i:=[\sigma_i]$ depends only on the restriction of $\omega_i$ to $\pi^{-1} (y_i)$ (and thus is independent of $\eta$).
This immediately yields that $b_i\cup \pi^* u =0$ for any $i$ and any $u\in H^* (N;\R)$, $\textrm{deg}\ u >0$, and that
\[
[\heta] = \pi^* [\eta] + \sum_{i=1}^m \delta_i b_i = \pi^* v + \sum_{i=1}^m \delta_i b_i\in \Kah (\hN,\hJ),
\]
which finishes the proof.
\endproof

\bigskip
\noindent
{\bf Proof of \ref{lem-gluing-two-psh-functions}.}

The form $\omega_i$ on $B(R)\setminus 0$ has locally integrable coefficients near $0$ and
thus defines a real degree-$(1,1)$ strictly positive current {\it on the whole $B(R)$} which is smooth on $B(R)\setminus 0$.

By the local $\partial\bar{\partial}$-lemma, we may assume, without loss of generality, that $\eta$ can be written on $B(R)$ as $\eta=\sqrt{-1}\partial\bar{\partial}F$ for a smooth strictly plurisubharmonic function $F$ on $B(R)$ and that
$T_{\omega_i} =\sqrt{-1}\partial\bar{\partial} T_G$ for some strictly plurisubharmonic function $G$ on $B(R)$ which is smooth
on $B(R)\setminus 0$. In particular, this means that $\omega_i=\sqrt{-1}\partial\bar{\partial}G$ on $B(R)\setminus 0$.

Define a function $H$ on $B(R)\setminus 0$ as $G$ if $G (0) = -\infty$ and as $\log |z|^2$ if $G (0)\in\R$. Recall that $\log |z|^2$ is a plurisubharmonic function of $z$ on the whole $\C^n$ which is, of course, smooth on $\C^n\setminus 0$.

Note that $a|z|^2+b$ is a plurisubharmonic function of $z$ for any $a,b\in \R$, $a>0$. Choose $a,b\in\R$, $a>0$, so that
\[
\min_{|z|=R/2} a|z|^2 +b > \max_{|z|=R/2} H, \ \max_{|z|=R/4} a|z|^2 +b < \min_{|z|=R/4} H.
\]
Then for a sufficiently small $\epsilon>0$ the regularized maximum $\max_\epsilon \{  a|z|^2 +b, H\}$ defines a strictly plurisubharmonic function on a neighborhood of the spherical annulus $\{ R/2\leq |z|\leq R/4\}$
which is equal to $a|z|^2 +b$ on a neighborhood of the sphere $\{ |z|=R/2\}$ and to $H$ on a neighborhood of the sphere $\{ |z|=R/4\}$.
Extending this function outside the sphere $\{ |z|=R/2\}$ by $a|z|^2 +b$ and inside the sphere $\{ |z|=R/4\}$ by $H$ we get a smooth strictly plurisubharmonic function $K$ on $B(R)\setminus 0$ such that $K (z) = a|z|^2 +b$ outside the ball $B(R/2)$ and $K = H$ on $B(R/4)\setminus 0$.

Since $F$ is strictly plurisubharmonic on $B(R)$ there exists a small $\varepsilon>0$ such that $F-\varepsilon (a|z|^2+b)$ is strictly plurisubharmonic on $B(3R/4)$. Thus, $L:=F-\varepsilon (a|z|^2+b) +\varepsilon K$ is a smooth function on $B(R)\setminus 0$ which, being a sum of two plurisubharmonic functions on $B(3R/4)\setminus 0$, is strictly plurisubharmonic on $B(3R/4)\setminus 0$. Note that $L=F$ outside $B(3R/4)$ and thus is strictly plurisubharmonic also there. Thus, $L$ is a smooth strictly plurisubharmonic function on $B(R)\setminus 0$ equal to $F$ outside $B(3R/4)$.

Observe that, since $K=H$ on $B(R/4)\setminus 0$ and $F-\varepsilon (a|z|^2+b)$ is continuous on $B(R)$, there exists $C_1>0$ such that
\begin{equation}
\label{eqn-mu-near-zero}
\varepsilon H - C_1\leq  L \leq \varepsilon H + C_1\ \textrm{on}\ B(R/4)\setminus 0.
\end{equation}

If $G (0) = -\infty$ and, accordingly, $H=G$, then, in view of \eqref{eqn-mu-near-zero} and since $G$ is continuous on $B(R)\setminus 0$, for any $0<\delta_i<\varepsilon$ one can find $C_2>0$ so that $\delta_i H - C_2=\delta_i G - C_2 < L$ outside $B(3R/4)$ and $\delta_i G - C_2 < L$ on $B(r)\setminus 0$ for some $0<r<R/4$. Therefore for a sufficiently small $\epsilon>0$ the regularized maximum $\varrho:=\max_\epsilon \{ L, \delta_i G - C_2\}$ is a smooth strictly plurisubharmonic function on $B(R)\setminus 0$ equal to $L$, and hence to $F$, outside $B(3R/4)$ and to $\delta_i G - C_2$ near $0$.

If $G (0)\in\R$ (meaning that $G$ is bounded from below on a neighborhood of $0$) and, accordingly, $H=\log |z|^2$, then, since $G$ is continuous on $B(R)\setminus 0$, for any $\delta_i >0$ one can find $C_3>0$ so that
$\delta_i G - C_3 < F=L$ outside $B(3R/4)$. Since $G$ is bounded from below
and $L=H=\log |z|^2$ near $0$ we get that $\delta_i G - C_3 > L$ on some neighborhood of $0$.
Therefore for a sufficiently small $\epsilon>0$ the regularized maximum $\varrho:=\max_\epsilon \{ L, \delta_i G - C_3\}$ is a smooth strictly plurisubharmonic function on $B(R)\setminus 0$ equal to $L$, and hence to $F$, outside $B(3R/4)$ and to $\delta_i G - C_3$ near $0$.

Let us sum up: in both cases ($G (0) = -\infty$ and $G (0)\in\R$), for any sufficiently small $\delta_i>0$ we get a smooth strictly
plurisubharmonic function $\varrho$ on $B(R)\setminus 0$ equal to $F$ outside $B(3R/4)$ and to $\delta_i G - C$, for some $C\in \R$, near $0$. Therefore for any sufficiently small $\delta_i>0$ we get a K\"ahler form $\xi_{i,\delta_i} :=\sqrt{-1}\partial\bar{\partial}\varrho$
on
$B(R)\setminus 0$ which equals $\delta_i \omega_i$ near $0$ and $\eta$ outside $B(3R/4)$.
\endproof

%%%%%%%%%%%%%%%%%%%%%%%%%%%%%%%%%%%%%%%%%%%%%%%%%%%%%%%%%%%%%%%%%%%%%%%%
\subsection{Weighted blow-ups and symplectic embeddings of ellipsoids}
\label{_weighted-blowup_Section}
%%%%%%%%%%%%%%%%%%%%%%%%%%%%%%%%%%%%%%%%%%%%%%%%%%%%%%%%%%%%%%%%%%%%%%%%

We are going to consider weighted blow-ups of smooth manifolds at a point. This operation can be performed both in complex and symplectic categories. Since we are going to compare weighted blow-ups of the
same smooth manifold with different complex structures, we will adapt the following point of view.

Fix a complex manifold $M$, $\dim_\C M =: n>1$ and a pairwise coprime vector $\ba := (a_1,\ldots,a_n)$. Note that the vector
$(1,\ba)=(1,a_1,\ldots,a_n)$ is then also pairwise coprime.

Topologically, a weighted blow-up of $M$ with the weight $(1,\ba)$ is defined as a topological space $\tM$ obtained by taking the connected
sum of $M$ with the weighted
projective space $\C P^n (1,\ba):= \C P^n (1,a_1,\ldots,a_n)$ near a regular point of $\C P^n (1,\ba)$. Recalling
\ref{_weighted-proj-space_Example_}, one easily sees that $\tM$ can
be equipped with the structure of a real orbifold.

A weighted blow-up at a point $x\in M$ can be realized in the complex
(or algebraic\footnote{In the algebraic language of modern algebraic geometry (see e.g.
\cite{_Hartshorne_}) the weighted blow-up of $M$ at a point $x$ with the weights $(1,a_1,\ldots,a_n)$ can be described as follows:
For each $k=1,\ldots,n$ assign the weight $a_k$ to the $k$-th coordinate $x_k$ on a neighborhood of $x$ in $M$ and let $R^m$ be the
ideal of polynomials of $x_1,\ldots,x_k$ of weighted degree $\geq m$. Consider $\bigoplus_{m=0}^\infty R^m$ as a graded scheme relative to $M$.
The graded spectrum (``Proj'') of this sheaf of rings gives the weighted blow-up
together with a projective morphism to $M$.})
category similarly to the usual blow-up (see e.g.
\cite[Sec. 6.38]{_Kollar-et-al-Rational-and-nearly-rational-2004_}) --
the latter can be viewed as the weighted blow-up with the weight
$(1,\ldots,1)$. Accordingly, any complex structure $I$ on $M$
defines a complex structure $\tI$ on $\tM$ and an
$(\tI,I)$-holomorphic map $\Pi_I: \tM\to M$ so that over $M\setminus
x$ the map $\Pi_I$ is a bi-holomorphism, while $\Pi_I^{-1} (x)=:
E(I)$ -- {\bf the exceptional divisor defined by $I$} -- is a
complex suborbifold biholomorphic (as an orbifold) to $\C P^{n-1}
(\ba)$. The complex structure $\tI$, the map $\Pi_I$ and
the exceptional divisor $E(I)$ are defined uniquely, up to a smooth
isotopy. The singular locus of $\tM$ is exactly the singular locus
of $E(I)$, that is, a finite collection of points.

If $J$ is another complex structure on $M$, we get another complex structure $\tJ$ on $\tM$
with another projection $\Pi_J: \tM\to M$ which is smoothly orbifold isotopic to $\Pi_I$ and
therefore induces the same map on cohomology which is independent of the
complex structure and will be denoted by $\Pi^*$. The exceptional divisor
$E(J)$ defined by $J$ might be different from $E(I)$
but lies in the same homology class which is independent of the complex structure.
We will denote the cohomology class that is Poincar\'e-dual to this homology class
by $e\in H^2 (\tM;\Z)$.

Now let us briefly recall how the weighted blow-up can be realized in the
symplectic category -- for more details see e.g. \cite{_Godinho_}.

For a pairwise coprime $\ba:=(a_1,\ldots,a_n)$ and $r>0$, denote by
$\cE_\ba (r)$ the ellipsoid $\pi\sum_{i=1}^n a_i |z_i|^2 \leq r$ in
the standard symplectic $\C^n$ with the coordinates
$z_1,\ldots,z_n$.

Given a symplectic manifold $(M,\omega)$, $\dim_\R M = 2n$, and
a symplectic embedding $\iota: \cE_\ba (r) \to (M,\omega)$, one can
construct an orbifold, which is orbifold diffeomorphic to $\tM$, by
removing $\iota (\cE_\ba (r))$ from $M$ and contracting the boundary
of the resulting manifold along the fibers of the $S^1$-action
induced by $\iota$ from the $S^1$-action on $\partial \cE_\ba (r)$
given by \eqref{_eqn-S1-action_}. The form $\omega$ is then extended
in a {\it canonical way} from $M\setminus \partial (\iota (\cE_\ba
(r)))$ to a symplectic form $\tomega$ on the orbifold $\tM$ -- this
procedure is called a symplectic cut (see \cite{_Lerman-MRL1995_}, \cite{_Godinho_},
cf. \cite{_Niederk-Pasq-JSG2009_}).

A calculation similar to \ref{_weighted-proj-space_Example_} shows that the cohomology class of $\tomega$ is given by
\[
[\tomega] = \Pi^* [\omega] - \frac{r}{\langle\ba\rangle} e.
\]

%%%%%%%%%%%%%%%%%%%%%%%%%%%%%%%%%%%%%%%%%%%%%%%%%%%%%%%%%%

\hfill

%%%%%%%%%%%%%%%%%%%%%%%%%%%%%%%%%%%%%%%%%%%%%%%%%%%%%%%%%%%%%%%%%%%%%%%%

The classical paper of McDuff-Polterovich \cite{_McD-Polt_} relates symplectic embeddings of
a disjoint union of $k$ balls into a symplectic manifold to the symplectic/K\"ahler classes of the blow-up
of the manifold at $k$ points. The proof of McDuff and Polterovich's results is based only on various local and global versions of
Moser's theorem and on local normal forms which generalize in a straightforward way to orbifolds (see \ref{_Moser_Proposition_},
\ref{_Moser_local_versions_Remark_}). Together with the construction of weighted blow-ups above this yields the following
version of McDuff and Polterovich's theorem
relating symplectic embeddings of a disjoint union of $k$ ellipsoids to the symplectic/K\"ahler classes of a weighted
blow-up.

\hfill

%%%%%%%%%%%%%%%%%%%%%%%%%%%%%%%%%%%%%%%%%%%%%%%%%%%%%%%%%%%%

\proposition (cf. \cite[Proposition 2.1A, 2.1B, 2.1C, Corollary 2.1D]{_McD-Polt_})
\label{_Kahler_McD_Polt_ellipses_Proposition_}\\
Let $M$ be a closed connected manifold equipped with a K\"ahler form $\omega$, $\dim_\C M=n$.
Let $\tM$ be a weighted blow-up of $M$ at $k$ points with weights $\ba_1,\ldots,\ba_k$.
Denote by $\Pi_I: \tM\to M$ the corresponding projection and by $e_1,\ldots,e_k\in H^2 (\tM;\R)$ the cohomology classes Poincar\'e-dual
to the homology classes of the corresponding exceptional divisors.

\smallskip
\noindent
A. For any sufficiently small positive $r_i$, $i=1,\ldots,k$, the symplectic manifold $(M,\omega)$
admits a symplectic embedding of a disjoint union of the ellipsoids
$\cE_{\ba_i} (r_i)$, $i=1,\ldots,k$, and for some complex structure $I$ on $M$ compatible with $\omega$
the cohomology class
\[
[\tomega] = \Pi^* [\omega] - \sum_{i=1}^k \frac{r_i}{\langle\ba_i\rangle} e_i\in H^2 (\tM;\R)
\]
is K\"ahler with respect to $\tI$.

\smallskip
\noindent
B. Assume there exists a complex structure $I$ of K\"ahler type on $M$
tamed by $\omega$ and a symplectic form $\tomega$ on $\tM$ taming
$\tI$ so that
\[
[\tomega] = \Pi^* [\omega] - \sum_{i=1}^k \frac{r_i}{\langle\ba_i\rangle} e_i
\]
for some $r_1,\ldots,r_k>0$.

Then $(M,\omega)$ admits a symplectic embedding of a disjoint union of the ellipsoids
$\cE_{\ba_i} (r_i)$, $i=1,\ldots,k$. \endproof

%%%%%%%%%%%%%%%%%%%%%%%%%%%%%%%%%%%%%%%%%%%%%%%%%%%%%%%%%%%%

\hfill

%%%%%%%%%%%%%%%%%%%%%%%%%%%%%%%%%%%%%%%%%%%%%%%%%%%%%%%%%%%%%

We will need the following version of \ref{_Kahler_McD_Polt_ellipses_Proposition_}, part B (cf. \cite[Theorem 8.3]{_EV-balls_}).

%%%%%%%%%%%%%%%%%%%%%%%%%%%%%%%%%%%%%%%%%%%%%%%%%%%%%%%%%%%%

\hfill

%%%%%%%%%%%%%%%%%%%%%%%%%%%%%%%%%%%%%%%%%%%%%%%%

\proposition\label{_McD_Polt_for_nearby_cs_Proposition_}\\
Let $(M, I, \omega)$, $\dim_\R M = 2n$, be a closed connected K\"ahler manifold. Let $k\in\N$ and let $\tM$,
$\Pi^* : H^2 (M;\R)\to H^2 (\tM;\R)$,  $[E_1], \ldots, [E_k]\in H^2(\tM, \Z)$,
$r_1,\ldots, r_k >0$ be as above.
Assume that there exists a complex structure $J$ of K\"ahler type on $M$ which is tamed by $\omega$ so that
\[
[\tomega] = \Pi^* [\omega] - \sum_{i=1}^k \frac{r_i}{\langle\ba_i\rangle} e_i\in \Kah (\tM, \tJ).
\]

Then $(M,\omega)$ admits a symplectic embedding of a disjoint union of the ellipsoids
$\cE_{\ba_i} (r_i)$, $i=1,\ldots,k$.

%%%%%%%%%%%%%%%%%%%%%%%%%%%%%%%%%%%%%%%%%%%%%%%%

\hfill

%%%%%%%%%%%%%%%%%%%%%%%%%%%%%%%%%%%%%%%%%%%%%%%%

\noindent {\bf Proof:}\\
The proof virtually repeats the proof of \cite[Theorem 8.3]{_EV-balls_}. We recall it briefly.

Note that $H^2 (\tM;\R) = H^2 (M;\R) \oplus \Span_\R \{ [E_1],\ldots,[E_k]\}$
and that the homomorphism $\Pi^* : H^2 (M;\R)\to H^2 (\tM;\R)$ acts as an identification of $H^2 (M;\R)$ with the first
summand. The classes $[E_1],\ldots,[E_k]$ are all of Hodge type $(1,1)$.
The identification preserves the Hodge types.

Since, by our assumption, $\Pi^*[\omega]_J^{1,1} - \sum_{i=1}^k \frac{r_i}{\langle\ba_i\rangle} e_i\in\Kah (\tM, \tJ)$,
it can be represented by a K\"ahler form $\talpha$ on $(\tM, \tJ)$.

Note that $\Pi^* [\omega]_J^{1,1}\in H^2 (\tM;\R)$ is of type $(1,1)$.
Hence, the class $\Pi^* [\omega] - \Pi^* [\omega]_J^{1,1}\in H^2 (\tM;\R)$ is of type $(2,0)+(0,2)$
and can be represented as $\Pi^* b$ for a $(2,0)+(0,2)$-class $b\in H^2 (M;\R)$. Represent $b$ by a
closed real-valued form $\beta$ on $M$ of type $(2,0)+(0,2)$ with respect to $J$. Then
the class $\Pi^* [\omega] - \Pi^* [\omega]_J^{1,1}$ is represented by a closed real-valued form
$\Pi_J^* \beta$ on $\tM$ of type $(2,0)+(0,2)$ with respect to $\tJ$.

Set $\tomega: = \talpha + \Pi_J^* \beta$. The form $\tomega$ is symplectic and tames $\tJ$.
The cohomology class of $\tomega$ can be written as
\[
[\tomega] = [\talpha] + [\Pi_J^* \beta] = \Pi^* [\omega]_J^{1,1} - \sum_{i=1}^k \frac{r_i}{\langle\ba_i\rangle} e_i + \Pi^* [\omega] - \Pi^* [\omega]_J^{1,1} =
\]
\[
=  \Pi^* [\omega] - \sum_{i=1}^k \frac{r_i}{\langle\ba_i\rangle} e_i.
\]
Now we can apply \ref{_Kahler_McD_Polt_ellipses_Proposition_} {\it with $J$ instead of $I$}, which yields the needed claim.
\endproof

%%%%%%%%%%%%%%%%%%%%%%%%%%%%%%%%%%%%%%%%%%%%%%%%

\hfill

%%%%%%%%%%%%%%%%%%%%%%%%%%%%%%%%%%%%%%%%%%%%%%%%%%%%%%%%%%%%
\proposition\label{_highest power_Proposition_}\\
With the notation as in \ref{_Kahler_McD_Polt_ellipses_Proposition_},
\[
\bigg\langle \bigg(\Pi^* [\omega] - \sum_{i=1}^k \frac{r_i}{\langle\ba_i\rangle} e_i\bigg)^n, [\tM]\bigg\rangle = \int_M \omega^n -  \sum_{i=1}^k \frac{r_i^n}{\langle\ba_i\rangle^{n}} =
\]
\[
= \int_M \omega^n - \Vol \bigg(\bigsqcup\limits_{i=1}^k \cE_{\ba_i} (r_i)\bigg).
\]

%%%%%%%%%%%%%%%%%%%%%%%%%%%%%%%%%%%%%%%%%%%%%%%%

\hfill

%%%%%%%%%%%%%%%%%%%%%%%%%%%%%%%%%%%%%%%%%%%%%%%%

\noindent {\bf Proof:}\\
Note that for all $i=1,\ldots,k$ we have $(\Pi^*
[\omega]) \cup e_i = 0$, as well as $e_i^n = -\langle\ba_i\rangle^{n-1}$, if $n$ is even, and
$e_i^n = \langle\ba_i\rangle^{n-1}$, if $n$ is odd (see \ref{_weighted-proj-space_Example_}). Also note that $e_i\cup e_j = 0$
for all $i\neq j$. Finally, recall that the symplectic volume of
$\cE_{\ba} (r)$ equals $r^{2n}/\langle\ba\rangle^n$. The claim follows directly from
these observations.
\endproof

%%%%%%%%%%%%%%%%%%%%%%%%%%%%%%%%%%%%%%%%%%%%%%%%

\hfill

%%%%%%%%%%%%%%%%%%%%%%%%%%%%%%%%%%%%%%%%%%%%%%%%%%%%%%%%%%%%
\section{Demailly-Paun theorem and the K\"ahler cone}
\label{_Demailly_Paun_Section_}
%%%%%%%%%%%%%%%%%%%%%%%%%%%%%%%%%%%%%%%%%%%%%%%%%%%%%%%%%%%%

%%%%%%%%%%%%%%%%%%%%%%%%%%%%%%%%%%%%%%%%%%%%%%%%

\hfill

%%%%%%%%%%%%%%%%%%%%%%%%%%%%%%%%%%%%%%%%%%%%%%%%
\proposition \label{_Kah_cone_blow_up_simple_Proposition_}\\
Assume $(M,I)$ is a closed connected Campana simple complex manifold, $\dim_\C M=n$.
Let $(\tM, \tI)$ be a weighted blow-up of $M$ at $k$ Campana-generic points $x_1,\ldots, x_k$ with weights $\ba_1,\ldots,\ba_k$ all of which are pairwise coprime. Define $\Pi_I: \tM\to M$, $E_i := E_i (I)=\Pi_I^{-1} (x_i)$ and $e_i\in H^2 (\tM;\Z)$, $i=1,\ldots,k$,
as above.

Assume that $\alpha\in \Kah (M,I)$.

Then, given $c_1,\ldots,c_k\in \R$,
the following claims are equivalent:
\begin{description}
\item[(A)] The cohomology class $\talpha:=\Pi^* \alpha - \sum_{i=1}^k c_i e_i \in H^2 (\tM;\R)$ is K\"ahler.
\hfill
\item[(B)] The conditions (B1) and (B2) below are satisfied:\\
(B1) All $c_i$ are positive. \\
(B2) $\langle \talpha^n, [\tM]\rangle >0$.
\end{description}

%%%%%%%%%%%%%%%%%%%%%%%%%%%%%%%%%%%%%%%%%%%%%%%%

\hfill

%%%%%%%%%%%%%%%%%%%%%%%%%%%%%%%%%%%%%%%%%%%%%%%%

Recall that $\tM$ is, in general, a non-smooth orbifold with isolated singularities.
In the case when $\tM$ is smooth (that is, in the case of usual, not weighted, blowups)
a claim similar to \ref{_Kah_cone_blow_up_simple_Proposition_}
was proved in \cite{_EV-balls_} using the Demailly-Paun theorem that describes the K\"ahler cone of a closed K\"ahler manifold:

%%%%%%%%%%%%%%%%%%%%%%%%%%%%%%%%%%%%%%%%%%%%%%%%

\hfill

%%%%%%%%%%%%%%%%%%%%%%%%%%%%%%%%%%%%%%%%%%%%%%%%
\theorem (Demailly-Paun, \cite{_Demailly_Paun_})
\label{_Dem_Paun_cone_Theorem_}\\
Let $X$ be a closed connected K\"ahler manifold. Let $\cK(X)\subset H^{1,1}(X;\R)$ be
the subset consisting of all (1,1)-classes $\zeta$ which satisfy
$\langle \zeta^s, [Z]\rangle >0$ for any homology class $[Z]$ realized by a complex subvariety
$Z\subset X$ of complex dimension $s>0$. Then the K\"ahler cone of $X$ is one of the
connected components of $\cK(X)$. \endproof

%%%%%%%%%%%%%%%%%%%%%%%%%%%%%%%%%%%%%%%%%%%%%%%%

\hfill

%%%%%%%%%%%%%%%%%%%%%%%%%%%%%%%%%%%%%%%%%%%%%%%%

\ref{_Dem_Paun_cone_Theorem_} cannot be directly applied to orbifolds (its orbifold version
seems to be true but its proof is not published, as far as we can ascertain) and therefore in order to prove
\ref{_Kah_cone_blow_up_simple_Proposition_}
we use an indirect argument where
\ref{_Dem_Paun_cone_Theorem_} is applied not to $\tM$ but to its K\"ahler resolution
constructed in part B of \ref{_pushforward_of_Kahler_classes_resolution_of_orbifolds_Theorem_}.

%%%%%%%%%%%%%%%%%%%%%%%%%%%%%%%%%%%%%%%%%%%%%%%%

\hfill

%%%%%%%%%%%%%%%%%%%%%%%%%%%%%%%%%%%%%%%%%%%%%%%%

\noindent
{\bf  Proof of \ref{_Kah_cone_blow_up_simple_Proposition_}.}\\

\noindent
{\bf  Proof of (A) $\Rightarrow$ (B).}\\
The implication (A) $\Rightarrow$ (B2) is obvious. To prove (A) $\Rightarrow$ (B1) note that, since $\talpha$ is K\"ahler,
for each $i=1,\ldots,k$ we have
$$0<\int_{E_i} \talpha^{n-1} = \int_{E_i} (-c_i e_i)^{n-1},$$
and since the restriction of $-e_i$ to $E_i$ is a positive multiple of the cohomology class of the
restriction of the Fubini-Study form $\Omega_{\ba_i}$ to the exceptional divisor ${E_i}$
and the integral
of the exterior power of the latter form over $E_i$ is positive, we readily get
 that $c_i >0$.
\endproof

%%%%%%%%%%%%%%%%%%%%%%%%%%%%%%%%%%%%%%%%%%%%%%%%

\hfill

%%%%%%%%%%%%%%%%%%%%%%%%%%%%%%%%%%%%%%%%%%%%%%%%

\noindent {\bf Proof of (B) $\Rightarrow$ (A).}\\
Assume (B1) and (B2) are satisfied.

Note that $\tM$ is a closed complex orbifold with isolated singularities and, by part A of \ref{_Kahler_McD_Polt_ellipses_Proposition_},
the complex structure $\tI$ on $\tM$ is of K\"ahler type. Denote by $y_1,\ldots, y_m\in \tM$ the singular points of $\tM$ -- each of them lies
in some exceptional divisor $E_i$. For each $i=1,\ldots,k$ denote by $\cS_i$ the set of $j$ such that $y_j\in E_i$.

Applying part B of \ref{_pushforward_of_Kahler_classes_resolution_of_orbifolds_Theorem_}
to $N:=\tM$ we get a K\"ahler resolution
$\pi: (\hN,\hI)\to (\tM,\tI)$ of $N=\tM$ and the cohomology classes $b_1,\ldots, b_m\in H^2 (\hN;\R)$ corresponding to $y_1,\ldots,y_m$.

Consider the map $\hpi:= \Pi_I\circ\pi : (\hN,\hI)\to (M,I)$. It is a biholomorphism over $M\setminus \{ x_1,\ldots,x_k\}$.

Consider the following family of cohomology classes of $\hN$:
\begin{equation}
\label{eqn-heta-lambda}
\hpi^* [\omega] - \lambda\sum_{i=1}^k c_i \pi^* e_i + \delta (b_1 +\ldots + b_m) = [\beta_{\lambda,\delta}]\in H^2 (\hN;\R),
\end{equation}
where $\lambda, \delta\geq 0$ and $\beta_{\lambda,\delta}$ is a smooth closed 2-form on $\hN$.
By part A of \ref{_Kahler_McD_Polt_ellipses_Proposition_}, $\Pi^* [\omega] - \lambda\sum_{i=1}^k c_i e_i\in \Kah (\tM,\tI)$ for any sufficiently small $\lambda >0$ and therefore, by part B of \ref{_pushforward_of_Kahler_classes_resolution_of_orbifolds_Theorem_},
for any sufficiently small $\lambda, \delta > 0$
the class $[\beta_{\lambda,\delta}]$ is K\"ahler for $(\hN, \hI)$ and the form $\beta_{\lambda,\delta}$ can be assumed to be K\"ahler.

Consider the cohomology class $\beta_{1,\delta} = \pi^* \talpha + \delta (b_1 +\ldots + b_m)$.
We claim that $\beta_{1,\delta} \in \cK (\hN)$ for sufficiently small $\delta>0$.

Indeed, let $Z$, $\dim Z =: s >0$, be a complex subvariety  of $(\hN, \hI)$.

There are 2 cases to consider: $Z\subsetneq \hN$ (Case I) and $Z=\hN$ (Case II).

Since $x_i$ are Campana-generic,
any connected
proper complex subvariety $Z$, $\dim Z =: s >0$, of $(\hN, \hI)$ either does not intersect any of the sets $\hpi^{-1} (x_i)$ (Case Ia) or is contained
in $\hpi^{-1} (x_i)$ for some $i$ (Case Ib).

\smallskip
\noindent
{\bf Case Ia:}

If $Z$ does not intersect any of the sets $\hpi^{-1} (x_i)$, then
\[ \langle \beta_{1,\delta}^s, [Z]\rangle = \langle \alpha^s, Z\rangle >0
\]  for all $\delta$, since $\alpha\in \Kah (M,I)$. This finishes the verification of the claim in the case (Ia).

\smallskip
\noindent
{\bf Case Ib:}

Assume $Z\subset \hpi^{-1} (x_i)$ for some $i$. Note that
\[
\hpi^{-1} (x_i) = \pi^{-1} (E_i) = \pi^{-1} \big( E_i \setminus \cup_{j\in\cS_i} \{ y_j\})\big)\cup \bigcup_{j\in\cS_i} \pi^{-1} (y_j).
\]
Also note that $\hpi^* \alpha$, all $e_j$, $j\neq i$, and all $b_j$ that correspond to $y_j\notin E_i$ vanish on $\hpi^{-1} (x_i)$. Hence
\begin{equation}
\label{eqn-talpha-on-Z-inside-a-fiber}
\langle \beta_{1,\delta}^s, [Z]\rangle = \bigg\langle \bigg( c_i \pi^* (-e_i) + \delta \sum_{j\in\cS_i} b_j\bigg)^s  , Z\bigg\rangle.
\end{equation}
The restriction of the class $-e_i$ to $E_i$ can be represented by the form $\Omega_i|_{E_i}$, where the $\Omega_i$ is the Fubini-Study orbifold K\"ahler form on $\C P^n (1,\ba_i)$ (recall that $\tM$ is obtained as a connected sum of $M$ with the weighted projective spaces $\C P^n (1,\ba_i)$, $i=1,\ldots,k$). Accordingly,
the restriction of the class $\pi^* (-e_i)$ to $\pi^{-1} (E_i)$ can be represented by the smooth form $\pi^* (\Omega_i|_{E_i})$.
The latter form is K\"ahler outside $\pi^{-1} (\cup_{j\in\cS_i} \{ y_j\})$
and its restriction to ${\pi^{-1} (y_j)}$ is zero for any $y_j\in E_i$.
Thus, for any $i$ we have $\pi^* (\Omega_i|_{E_i})\geq 0$ on $\pi^{-1} (E_i)$, hence
\[
\langle (-e_i)^s, [Z]\rangle \geq 0.
\]
Recall that for any sufficiently small $\lambda, \delta >0$ the form $\beta_{\lambda,\delta}$ is K\"ahler on $(\hN, \hI)$.
Therefore for any sufficiently small $\lambda, \delta >0$ (independent of $Z$!) the form $c_i\pi^* (\Omega_i|_{E_i}) + \beta_{\lambda,\delta}$ is positive on $\pi^{-1} (E_i)$ for any $i$ and therefore so is its restriction to the complex subvariety $Z\subset \pi^{-1} (E_i)$.
Hence,
\[
0< \langle [c_i\pi^* (\Omega_i|_{E_i}) + \beta_{\lambda,\delta}]^s, [Z]\rangle.
\]
On the other hand, in view of \eqref{eqn-heta-lambda}, the cohomology class of the restriction of $c_i\pi^* (\Omega_i|_{E_i}) + \beta_{\lambda,\delta}$ to $\pi^{-1} (E_i)$
can be written as
\[
[c_i\pi^* (\Omega_i|_{E_i}) + \beta_{\lambda,\delta}] = (1+\lambda) c_i \pi^* (-e_i) + \delta \sum_{j\in\cS_i} b_j.
\]
Thus, for any sufficiently small $\lambda, \delta >0$ (independent of $Z$)
\[
0< \bigg\langle \bigg( (1+\lambda) c_i \pi^* (-e_i) + \delta \sum_{j\in\cS_i} b_j\bigg)^s, [Z]\bigg\rangle.
\]
Note that, by the properties of the classes $b_j$ given by part B of \ref{_pushforward_of_Kahler_classes_resolution_of_orbifolds_Theorem_}, $\pi^* (-e_i) b_j = 0$ for any $i,j$. Therefore for any sufficiently small $\lambda, \delta >0$ (independent of $Z$)
\[
0< \bigg\langle \bigg( (1+\lambda) c_i \pi^* (-e_i) + \delta \sum_{j\in\cS_i} b_j\bigg)^s, [Z]\bigg\rangle =
\]
\[
\bigg\langle \bigg( c_i \pi^* (-e_i) + \delta \sum_{j\in\cS_i} b_j \bigg)^s, [Z]\bigg\rangle + \bigg(\sum_{j=1}^s \binom{s}{j}\lambda^j\bigg) \bigg\langle (-e_i)^s, [Z]\bigg\rangle.
\]
Since $\lambda^j \langle (-e_i)^s, [Z]\rangle\geq 0$ for any $\lambda >0$ and any $j\in \N$, we get, by \eqref{eqn-talpha-on-Z-inside-a-fiber}, that for any sufficiently small $\delta >0$ (independent of $Z$)
\[
\langle \beta_{1,\delta}^s, [Z]\rangle = \bigg\langle \bigg(c_i \pi^* (-e_i) + \delta \sum_{j\in\cS_i} b_j\bigg)^s  , Z\bigg\rangle >0.
\]
This finishes the verification of the claim in the case (Ib).

\smallskip
\noindent
{\bf Case II:}

Assume $Z=\hN$. Then, by (B2) and since $\pi$ is of degree $1$,
\[
\langle (\pi^*\talpha)^n, [\hN]\rangle = \langle \talpha^n, [\tM]\rangle >0.
\]
Therefore, since $\beta_{1,\delta}\to \pi^*\talpha$ as $\delta\to 0$, for any sufficiently small $\delta >0$
\[
\langle \beta_{1,\delta}^n, [\hN]\rangle >0.
\]
This finishes the verification of the claim in the case (II).

\smallskip
Thus we have proved that $\beta_{1,\delta} \in \cK (\hN)$ for a sufficiently small $\delta>0$.

Let us now show that there exists a K\"ahler class
in the connected component of $\cK(\hN)$ containing $\beta_{1,\delta}$.

Indeed, similarly to \ref{_highest power_Proposition_}, one gets that (B2) is equivalent to the condition
\[
\sum_{i=1}^k c_i^n < \langle \alpha^n, [M]\rangle,
\]
If this condition holds for $c_1,\ldots, c_k >0$, it also holds for
$\epsilon c_1,\ldots, \epsilon c_k$ for any $\epsilon\in (0,1]$. For any such $\epsilon$ the
numbers $\epsilon c_1,\ldots, \epsilon c_k$ are still positive and
therefore, by the argument above, for any sufficiently small $\delta >0$ and any $\epsilon\in (0,1]$ the
class $\hpi^* \alpha - \epsilon \sum_{i=1}^k c_i \pi^* e_i + \epsilon \delta (b_1+\ldots+b_m)$ also lies in $\cK(\hN)$. But, as we have already
mentioned above, for any sufficiently small $\epsilon >0$
the class
\[
\hpi^* \alpha - \epsilon \sum_{i=1}^k c_i \pi^* e_i + \epsilon \delta (b_1+\ldots+b_m) = [\beta_{\epsilon, \epsilon\delta}]
\]
is K\"ahler.
Thus, for any sufficiently small $\delta>0$ the class
$\beta_{1,\delta}$ lies in the same connected component of $\cK(\hN)$ as a
K\"ahler class $[\beta_{\epsilon, \epsilon\delta}]$. Therefore, by
\ref{_Dem_Paun_cone_Theorem_}, the class
\[
\beta_{1,\delta} = \pi^* \talpha + \delta (b_1 +\ldots + b_m).
\]
is K\"ahler on $(\hN,\hI)$.
Since $\pi$ is a biholomorphism outside $\pi^{-1} \{ y_1,\ldots,y_m\}$ and $\pi_* b_j =0$ for any $j$, by part A of \ref{_pushforward_of_Kahler_classes_resolution_of_orbifolds_Theorem_},
$\pi_* \beta_{1,\delta} = \talpha$
is a K\"ahler class on $(\tM,\tI)$, which finishes the proof.
\endproof

%%%%%%%%%%%%%%%%%%%%%%%%%%%%%%%%%%%%%%%%%%%%%%%%

\hfill

%%%%%%%%%%%%%%%%%%%%%%%%%%%%%%%%%%%%%%%%%%%%%%%%

\section{Proof of \ref{_USP_Campana_Simple_Theorem_}}
\label{_USP_Campana_Simple_Theorem_Pf_Section_}

Let us say that a closed ellipsoid is {\bf simple} if it is of the form $\cE_\ba (r)$ for a pairwise coprime vector $\ba$ and some $r>0$.

Any vector with positive coordinates can be approximated by vectors proportional to pairwise coprime vectors (see \ref{_approximation-by-pairwise-coprime-vectors_Proposition_} below) and
therefore any open neighborhood of any closed ellipsoid contains a simple closed ellipsoid. Thus it suffices to prove that if there is a collection of disjoint simple closed ellipsoids $\cE_{\ba_i} (r_i)$, $i=1,\ldots,k$, whose total volume is less then $\Vol M$, then there exists a symplectic embedding of their union into $(M,\omega)$.

Consider a disjoint union $\bigsqcup\limits_{i=1}^k \cE_{\ba_i} (r_i)$ whose total symplectic volume is less than the symplectic volume of $M$, that is,
\begin{equation}
\label{eqn-total-vol-balls-less-than-vol-M}
\sum_{i=1}^k r_i^{n}/\langle\ba_i\rangle^n < \Vol (M) = \langle [\omega]^n, [M]\rangle.
\end{equation}
We need to show that it admits a symplectic embedding into $(M,\omega)$.

It follows from the Kodaira-Spencer stability theorem \cite{_Kod-Spen-AnnMath-1960_} (see \cite[Theorem 5.6]{_EV-balls_} for more details) and the hypothesis of the theorem
that there exists a
Campana simple complex structure
$J$ on $M$ sufficiently close to $I$ with the following properties:

\smallskip
\noindent
1. $[\omega]_J^{1,1}\in \Kah (M,J)$ (this follows from the Kodaira-Spencer stability theorem -- see \cite[Theorem 5.6]{_EV-balls_}),

\smallskip
\noindent
2.
\begin{equation}
\label{eqn-total-vol-balls-less-than-vol-M-J}
\sum_{i=1}^k r_i^{n}/\langle\ba_i\rangle^n  < \big\langle ([\omega]_J^{1,1})^n, [M]\big\rangle.
\end{equation}
(This is possible by \eqref{eqn-total-vol-balls-less-than-vol-M}).

\smallskip
\noindent
3. $J$ is tamed by $\omega$ (this is possible because $I$ is tamed by $\omega$),

\smallskip
\noindent
Choose $k$ distinct Campana-generic points $x_1,\ldots,x_k \in (M,J)$,
and consider the weighted blow-up $(\tM, \tJ)$
of $(M,J)$ at those points with the weights $\ba_1,\ldots,\ba_k$. By
\ref{_Kah_cone_blow_up_simple_Proposition_} applied to the K\"ahler
class $[\omega]_J^{1,1}$, the cohomology class $\Pi^*
[\omega]_J^{1,1} - \sum_{i=1}^k r_i/\langle\ba_i\rangle e_i$ is K\"ahler with respect
to $\tJ$ (note that, by \ref{_highest power_Proposition_}, the
condition (B2) in \ref{_Kah_cone_blow_up_simple_Proposition_} is
equivalent to \eqref{eqn-total-vol-balls-less-than-vol-M-J}).
Therefore, by \ref{_McD_Polt_for_nearby_cs_Proposition_}, $(M,\omega)$
admits a symplectic embedding of $\bigsqcup\limits_{i=1}^k \cE_{\ba_i} (r_i)$.
\endproof

%%%%%%%%%%%%%%%%%%%%%%%%%%%%%%%%%%%%%%%%%%%%%%%%

\hfill

%%%%%%%%%%%%%%%%%%%%%%%%%%%%%%%%%%%%%%%%%%%%%%%%

\section{Appendix: Approximation by pairwise coprime vectors}
\label{_approximation_appendix_Section_}

In this appendix we will prove that any vector with positive coordinates can be approximated by vectors proportional to
pairwise coprime vectors. The result is probably known but we have been unable to find it in the literature. In fact, below we
present a proof (due to Uri Shapira) of a stronger claim.

%%%%%%%%%%%%%%%%%%%%%%%%%%%%%%%%%%%%%%%%%%%%%%%%

\hfill

%%%%%%%%%%%%%%%%%%%%%%%%%%%%%%%%%%%%%%%%%%%%%%%%

\proposition
\label{_approximation-by-pairwise-coprime-vectors_Proposition_}\\
Any vector with positive coordinates in $\R^n$ can be approximated by vectors proportional to vectors whose coordinates are pairwise different primes.

%%%%%%%%%%%%%%%%%%%%%%%%%%%%%%%%%%%%%%%%%%%%%%%%

\hfill

%%%%%%%%%%%%%%%%%%%%%%%%%%%%%%%%%%%%%%%%%%%%%%%%

\noindent
{\bf Proof (Uri Shapira):}\\
Consider a vector $(x_1,\ldots,x_n)\in \R^n$ with positive coordinates. Take an arbitrary $\varepsilon >0$ and let us find a vector at the $l_\infty$-distance $\leq \varepsilon$ from $(x_1,\ldots,x_n)$ which is proportional to a vector whose coordinates are pairwise different primes.

Choose a vector $(y_1,\ldots,y_n)$
with positive pairwise different rational coordinates so that $\max_i |y_i-x_i|\leq \varepsilon/2$. Let
\[
c:=\min_{i,j, i\neq j} |y_i-y_j|>0,
\]
\[
C:= \max_i y_i >0
\]
According to a theorem of Hoheisel \cite{_Hoheisel_}, there exists $0<\vartheta <1$ and $l_0\in\N$ so that for any $l\geq l_0$ the interval $(l,l+l^\vartheta)$ contains a prime. Choose $l_1\geq l_0$ so that
\[
l^\vartheta < l\min \bigg\{ \frac{\varepsilon}{2C}, c\bigg\}\ \textrm{for any}\ l\geq l_1.
\]
Choose a sufficiently large $N\in \N$ so that $Ny_1,\ldots,Ny_n$ are integers greater than $l_1$. Then the intervals
$(Ny_i, Ny_i+ (Ny_i)^\vartheta)$, $i=1,\ldots,n$, are pairwise disjoint and each of them contains a prime: $p_i\in (Ny_i, Ny_i+ (Ny_i)^\vartheta)$. Since the intervals are disjoint, the primes $p_1,\ldots,p_n$ are pairwise different, and since $p_i\in (Ny_i, Ny_i+ (Ny_i)^\vartheta)$, we get
\[
\max_i |p_i/N - y_i| \leq \max_i \frac{(Ny_i)^\vartheta}{N} \leq \max_i \frac{(\varepsilon/2C) Ny_i}{N}\leq \varepsilon/2.
\]
Thus, the vector $(p_1/N,\ldots, p_n/N)$ is proportional to the vector $(p_1,\ldots,p_n)$ whose coordinates are pairwise different primes and
\[
\max_i |p_i/N - x_i|\leq \max_i |p_i/N - y_i| + \max_i |y_i - x_i| \leq \varepsilon/2 + \varepsilon/2 = \varepsilon,
\]
as required.
\endproof

%%%%%%%%%%%%%%%%%%%%%%%%%%%%%%%%%%%%%%%%%%%%%%%%%%%

\hfill

\noindent {\bf Acknowledgements:} We are grateful to L. Polterovich
for useful discussions and many interesting suggestions. We thank Uri Shapira for communicating to us the proof of \ref{_approximation-by-pairwise-coprime-vectors_Proposition_} and Felix Schlenk for enlightening comments on a preliminary version
of the paper.

\hfill

{\small
\noindent {\sc Michael Entov\\
Department of Mathematics,\\
Technion - Israel Institute of Technology,\\
Haifa 32000, Israel}\\
{\tt  entov@math.technion.ac.il}
}
\\

{\small
\noindent {\sc Misha Verbitsky\\
            {\sc Instituto Nacional de Matem\'atica Pura e
              Aplicada (IMPA) \\ Estrada Dona Castorina, 110\\
Jardim Bot\^anico, CEP 22460-320\\
Rio de Janeiro, RJ - Brasil}

\smallskip
\noindent
and

\smallskip
\noindent
Laboratory of Algebraic Geometry,\\
National Research University HSE, Faculty of Mathematics, \\
7 Vavilova Str., Moscow, Russia}\\
{\tt  verbit@mccme.ru}
}

\end{document}